\documentclass[a4paper,11pt,twoside]{article}

\usepackage[english]{babel}
\usepackage{amsfonts,amsmath,amssymb,amstext,amsthm}
\usepackage{lmodern,fancyhdr,lastpage,nccfoots,caption}
\usepackage{hyperref}
\usepackage[utf8]{inputenc}
\usepackage{graphicx,xcolor}
\usepackage{enumerate}
\usepackage[a4paper,margin=1in]{geometry}

\usepackage[all]{xy}
\xyoption{arc}
\CompileMatrices

\newtheorem{thm}{Theorem}[section]
\newtheorem{cor}[thm]{Corollary}

\newtheorem{prop}[thm]{Proposition}

\theoremstyle{definition}
\newtheorem{Def}[thm]{Definition}    
\newtheorem{ex}[thm]{Example}
\newtheorem{rem}[thm]{Remark}

\newcommand{\CC}{\mathbb{C}}
\newcommand{\RR}{\mathbb{R}}
\newcommand{\QQ}{\mathbb{Q}}

\newcommand{\sA}{\mathcal{A}}       
\newcommand{\sE}{\mathcal{E}}       
\newcommand{\gG}{\mathcal{G}}               
\newcommand{\sM}{\mathcal{M}}       
\newcommand{\sN}{\mathcal{N}}       
\newcommand{\sO}{\mathcal{O}}       

\newcommand{\lieg}{\mathfrak{g}}
\newcommand{\liegl}{\mathfrak{gl}}

\newcommand{\Gr}{\operatorname{Gr}}
\newcommand{\Ker}{\operatorname{Ker}}
\newcommand{\Lie}{\operatorname{Lie}}

\newcommand{\Hom}{\operatorname{Hom}}

\newcommand{\trace}{\operatorname{trace}}

\newcommand{\Symp}{\operatorname{Symp}}
\newcommand{\Diff}{\operatorname{Diff}}
\newcommand{\Vect}{\operatorname{Vect}}

\newcommand{\too}{\longrightarrow}

\newcommand{\SL}{\operatorname{SL}}

\DeclareMathOperator{\End}{End}     
\DeclareMathOperator{\GCD}{GCD}     
\DeclareMathOperator{\rk}{rk}       

\newcommand{\bS}{\mathbb{S}}        

\renewcommand{\geq}{\geqslant}        
\renewcommand{\leq}{\leqslant}        
\newcommand{\ox}{\otimes}             
\newcommand{\sse}{\subseteq}		
\renewcommand{\:}{\colon}             

\newcommand{\suchthat}{\;\;|\;\;}

\renewcommand{\leq}{\leqslant}
\renewcommand{\geq}{\geqslant}

\renewcommand{\phi}{\varphi}

\newcommand{\xra}{\xrightarrow}
\DeclareMathOperator{\HNF}{HNF}
\DeclareMathOperator{\HNT}{HNT}

\DeclareMathOperator{\Ad}{Ad}

\newcounter{a}
\ifodd\thea\else\stepcounter{a}\fi

\fancypagestyle{plain}{%
\fancyhf{}

}
\defineshorthand{"-}{\nobreak\hskip0pt\discretionary{-}{-}{-}\nobreak\hskip0pt}

\def\Circlearrowright{\ensuremath{%
  \rotatebox[origin=c]{180}{$\circlearrowright$}}}

\newcommand{\word}[1]{\quad\text{#1}\quad} 

\begin{document}
\thispagestyle{plain}

\begin{center}
\Large
\textsc{Geometric Invariant Theory,\\
holomorphic vector bundles and\\
the Harder--Narasimhan filtration}
\end{center}


\begin{center}
  {\em Alfonso Zamora}\\
  \small Departamento de Matem\'atica Aplicada y Estad\'istica\\
    \small Universidad CEU San Pablo\\
  \small Juli\'an Romea 23, 28003 Madrid, Spain\\
  \small e-mail: \texttt{alfonso.zamorasaiz@ceu.es}
\end{center}

\begin{center}
  {\em Ronald A. Z\'u\~niga-Rojas}\\
  \small Centro de Investigaciones Matem\'aticas y Metamatem\'aticas CIMM\\
  \small Escuela de Matem\'atica, Universidad de Costa Rica UCR\\
  \small San Jos\'e 11501, Costa Rica\\
  \small e-mail: \texttt{ronald.zunigarojas@ucr.ac.cr}
\end{center}

\textbf{Abstract.}  
This survey intends to present the basic notions of Geometric Invariant Theory (GIT) through its paradigmatic application in the construction of the moduli space of holomorphic vector bundles. Special attention is paid to the notion of stability from different points of view and to the concept of maximal unstability, represented by the Harder-Narasimhan filtration and, from which, correspondences with the GIT picture and results derived from stratifications on the moduli space are discussed. 

\medskip

\textbf{Keywords:} 
Geometric Invariant Theory, Harder-Narasimhan filtration, moduli spaces, vector bundles, Higgs bundles, GIT stability, symplectic stability, stratifications.

\textbf{MSC class:}	\texttt{14D07, 14D20, 14H10, 14H60, 53D30}

\tableofcontents

\section{Introduction}
\label{sec:intro}

Moduli spaces are structures classifying objects under some equivalence relation and  many of these problems can be posed as quotients of a projective variety $X$ under a reductive group $G$. The purpose of Geometric Invariant Theory (abbreviated
GIT, \cite{Mu, MFK}) is to provide a way to define a quotient of $X$ by the action of $G$ with an algebro-geometric structure.
This way, GIT results assure a good structure for the quotient, giving a positive solution to the classification problem.

Let $G$ be a reductive complex Lie group acting on an algebraic variety $X$. In the case when the variety $X$ is affine there is a simpler solution which dates back to Hilbert's $14^{th}$ problem.
Let $A(X)$ denote the coordinate ring of the affine variety $X$. Nagata \cite{Na} proved that if $G$ is reductive,
the ring of invariants $A(X)^{G}$ is finitely generated, hence is the coordinate ring of an affine variety, therefore we can
define the quotient of $X$ by $G$ as the affine variety associated to the ring $A(X)^{G}$. First thing to note is that the orbit space (i.e. the quotient space where each point corresponds to an orbit) is not separated, the reason why the ring of invariants identifies orbits under the notion $S$-equivalence (see Example \ref{hyperboles}), yielding a Hausdorff quotient. 

When taking the quotient of a projective variety $X$ by a group $G$, there are bigger issues which have to be taken into account.  Actions on the variety $X$ do not determine uniquely an action on its ring of functions, which makes necessary to pass through the affine cone $\hat{X}$ and the linearization of the action. The affine cone introduces an origin which will be removed when going back to the projective variety $X$, and those orbits falling through this origin will behave bad when taking the quotient by $G$. This is how the notion of GIT stability appears. 

GIT was developed by David Mumford \cite{Mu} as a theory to construct quotients of varieties, broadening the scope of Hilbert's problem. After Narasimhan and Seshadri celebrated theorem \cite{NS, Se} relating stable bundles to irreducible representations of the unitary group, GIT mayor application was the construction of a projective variety classifying all holomorphic structures in a smooth bundle, what is called the moduli space of vector bundles.
Based on Kirwan's work \cite{Ki}, there is another theory of symplectic quotients mirroring with the GIT picture, the Kempf-Ness theorem \cite{KN} being the link between GIT and symplectic stability.

Extending representations of the unitary group to the whole $G$ we find the notion of Higgs bundle, first studied by Hitchin \cite{Hi1} as the solutions of certain partial differential equations on a Riemann surface. These objects have turned to be a central element intertwining geometry, topology and physics after the works of Atiyah-Bott \cite{AB}, the generalization to higher dimension by Simpson \cite{Si2}, the study of its moduli space by Hausel \cite{Hau}, the moduli construction using GIT by Nitsure \cite{Ni} or its Hitchin-Kobayashi correspondence in \cite{G-PGM}. 

When constructing the moduli space of vector bundles of the moduli space of Higgs bundles, we impose a stability condition and declare certain objets to be stable, semistable or polystable, such that when encoding the problem in the GIT framework, these notions match. Therefore, constructing a GIT quotient of $S$-equivalence classes of orbits we provide a moduli space classifying $S$-equivalence classes of semistable bundles, each one containing a polystable representative, as desired. This stability concept and the moduli theory leaves outside of the quotient picture the non-semistable objects, the unstable ones. Unstable bundles are more complex: they carry bigger and different automorphism groups and, therefore, the action of $G$ on them yields different stabilizers, which causes anomalies in the quotient, the reason we have to remove them. 

Hilbert-Mumford criterion is the tool to check GIT stability through the idea of 1-parameter subgroups, which are 1 dimensional subvarieties of the group $G$ accounting for the various features of separating orbits. GIT unstable points are localized as those contradicting the criterion, analogous to unstable bundles being the objects not verifying the stability condition for bundles. The 
result of Kempf \cite{Ke} and the well known Harder-Narasimhan filtration \cite{HN} show the two sides of this maximal principle which are put in correspondence in \cite{Za3, GSZ1}, and which is linked to symplectic and differencial geometry in \cite{GRS}. For a modern treatment generalizing the structure of unstability to quotient stacks problems see Halpern-Leistner theory \cite{H-L}.

One of the main applications of the Harder-Narasimhan filtration is that it provides a way to stratify the moduli problem in terms of the numerical invariants of the object, in order to study its properties. For vector bundles, Shatz \cite{Sha} constructs a locally closed stratification indexed by the data of the Harder-Narasimhan filtration, and Hoskins and Kirwan give a structure of moduli space to each strata through maximal GIT unstability in \cite{HK}. Particularly, for Higgs bundles, 
these stratifications provide a useful way to study the geometry and topology of the moduli space, see \cite{B-B, Hau, Z-R1, GZ-R}.

This paper intends to introduce the main elements of Geometric Invariant Theory through various examples, particularly the moduli space of holomorphic vector bundles. Special attention is devoted to unstable objects to provide a survey of results of the authors on correspondences with the GIT picture and stratifications of the moduli space of Higgs bundles.

The survey is intended to be self-contained and is organized as follows. In Section \ref{sec:preliminaries}
we give the basic definitions that we will need about
Lie groups and algebras, algebraic varieties and vector bundles. After that, in Section \ref{sec:git}, we introduce Geometric Invariant Theory and the notion of GIT stability. Subsection \ref{ssec:sympl-stability} links this theory to symplectic geometry and Subsection \ref{ssec:examples-stability} shows different examples through which these notions are visualized, the construction of the projective space and the Grassmannian variety as GIT or symplectic quotients, and the study
of the classical problem (dating back to Hilbert) of classifying configurations of points on the projective line. Section \ref{sec:ms} goes over the GIT construction of the moduli space of holomorphic vector bundles and presents the Harder-Narasimhan filtration in Subsection \ref{ssec:HN}. Section \ref{sec:Higgs} is devoted to the moduli space of Higgs bundles, from their algebraic and analytical points of view. Section \ref{AZresults} includes results on correspondences between Harder-Narasimhan fitrations and maximally destabilizing 1-parameter subgroups for unstable objects in different moduli problems. Finally, Section \ref{sec:RZRresults} shows results on stratifications on the moduli space of vector bundles. 

\subsection*{Acknowledgements}

We wish to thank to a number of people for their knowledge, expertise and stimulating conversations towards completing this work: Luis \'Alvarez-C\'onsul, \'Oscar Garc\'ia-Prada, Tom\'as G\'omez, Peter Gothen, Carlos Florentino, Alessia Mandini, Peter Newstead, Andr\'e Oliveira, Milena Pabiniak, Ignacio Sols.

We also thank the support of Universidad de Costa Rica through Escuela de Matem\'atica, specifically through Centro de Investigaciones Matem\'aticas 
y Metamatem\'aticas (CIMM), and through Oficina de Asuntos Internacionales y de Cooperaci\'on Externa (OAICE) through Programa Acad\'emicos Visitantes (PAV), for their support and the invitation to give a course on March 2016, where most part of this material was covered. 
  
The first author is supported by project {\tt MTM2016-79400-P} of the Spanish government. The second author is supported by Universidad de Costa Rica through Escuela de Matem\'atica and CIMM with projects {\tt 820-B5-202} and {\tt 820-B8-224}.

\section{Preliminaries}
\label{sec:preliminaries}

\subsection[Lie groups]{Lie groups}
\label{ssec:lie-groups}

\begin{Def}
 A {\em Lie group} $G$ is a group which is also a differentiable manifold such that the product map  
 \[
 *\: G \times G \to G  
 \]
 \[
  (g,h)\mapsto g*h
 \] 
 and the inverse map
 \[
 (\cdot)^{-1}\: G\to G 
 \]
 \[
 g\mapsto g^{-1} 
 \]
 are differentiable. If $G$ is a complex manifold and the operations $*$ and $(\cdot)^{-1}$ are analytic, we say that $G$ is a {\em complex Lie group}.
\end{Def}

Basic examples of real Lie groups are the additive real group $(\RR,+)$ which is connected and 
the multiplicative real group $\RR^{*} = (\RR-\{0\},\cdot)$, being  disconencted. 
The multiplicative complex group $\CC^{*} = (\CC-\{0\},\cdot)$ is a comlpex Lie group and so, the unitary circle $(\bS^{1},\cdot)$ is a complex Lie subgroup of $\CC^*$. By complexifying real varieties into complex ones, we can complexify real groups $K$ to obtain complex groups $K_{\CC}$. For example, $\CC^{*}$ is the complexification of the unitary circle $\bS^{1}$. 

Special importance is given to matrix Lie groups, those which can be seen as subgroups of the square $n\times n$ matrices, $M^{n\times n}(k)$:
 \begin{enumerate}[(a)]
  \item 
  The general linear real group
  \[
   GL(n,\RR) =
   \{
   A 
   \in M^{n\times n}(\RR)\suchthat 
   \det(A) \neq 0
   \}
  \]
  of $n\times n$ invertible real matrices is a real Lie group, and its complexification,
the general linear complex group
  \[
   GL(n,\CC) =
   \{
   A 
   \in M^{n\times n}(\CC)\suchthat 
   \det(A) \neq 0
   \}
  \]
  of $n\times n$ invertible complex matrices, is a complex Lie group. The special linear group
  \[
   SL(n,\CC) =
   \{
   A 
   \in GL(n,\CC)\suchthat 
   \det(A) = 1
   \}
  \]
  is a complex Lie subgroup of $GL(n,\CC)$ and the unitary group
  \[
   U(n) =
   \{
   A 
   \in GL(n,\CC)\suchthat 
   A A^{*} = A^{*} A = I_n
   \}
  \]
  is a compact Lie subgroup of $GL(n,\CC)$ whose complexification is  $U(n)_{\CC}=GL(n,\CC)$. The unitary circle $\bS^1$ can be identified with $U(1)$.
  
  \item 
   The orthogonal real group
  \[
   O(n,\RR) =
   \{
   A 
   \in GL(n,\RR)\suchthat 
   A A^{t} = A^{t} A = I_n
   \}
  \]
  and the special orthogonal real group
  \[
   SO(n,\RR) =
   \{
   A 
   \in \sO(n)\suchthat 
   \det(A) = 1
   \}
  \]
  are real matrix Lie groups, whose complexifications are $O(n,\CC)$ and $SO(n,\CC)$. For example, the special orthogonal real group of order $2$
  \[
   SO(2,\RR) =
   \left\{
   A = 
   \left(
   \begin{array}{c r}
    \cos(\theta) & -\sin(\theta)\\
    \sin(\theta) & \cos(\theta)
   \end{array}
   \right)
   \in GL(2,\RR)\ |\ 
   \theta \in [0,2\pi]
   \right\}
  \]
  corresponds to the group of rotations in $\RR^2$, preserving orientation. This group is also a complex Lie group by its own, since there is a 
  diffeomorphism
  \[
   \bS^{1} \xra{\quad \cong \quad} SO(2).
  \]
  
  \item 
  The symplectic group 
  \[
   Sp(2n,\CC) =
   \left\{
   A \in
   M^{2n\times 2n}(\CC) 
   \suchthat 
   A^{t} J A = J
   \right\}
  \]
  where
  \[
   J =
   \left(
   \begin{array}{c r}
     0  & -I_n\\
    I_n &   0
   \end{array}
   \right)
   \in
   M^{2n\times 2n}(\CC) 
  \]
  is the group of matrices preserving the standard symplectic form in $\CC^{2n}$. 
 \end{enumerate}

\subsection[Lie algebras]{Lie algebras}
\label{ssec:lie-algebras}

\begin{Def}
 A {\em Lie algebra} is a $k$-vector space $\lieg$ joint with a non-associative, bilinear, alternating map, called \emph{Lie bracket},
 \[
  [\cdot,\cdot]\:
  \lieg \times \lieg \to \lieg
 \]
 \[
  (v,w)\mapsto [v,w]\; ,
 \]
satisfying Jacobi's identity
 \[
  [u,[v,w]] + [v,[w,u]] + [w,[u,v]] = 0\; ,
 \]
 for any $u,v,w\in \lieg$.
\end{Def}

Any vector space $V$ over a field $k$ with the trivial bracket $[v,w] = 0$ is a Lie algebra over $k$. A non trivial example of a real Lie algebra is $\lieg = \RR^3$ with the Lie bracket $[v,w] := v\times w$.

\begin{ex}
Another example comes from an associative algebra $(A, *)$ over a field $k$ with the Lie bracket 
 \[
  [v,w] = v*w - w*v
 \]
 for any $v,w\in A$. In particular, if $V$ is an $n$-dimensional vector space over $k$, then $\lieg = \End_{k}(V)$ is an associative algebra with composition, 
 so it is a Lie algebra with the bracket $[f,g] = f\circ g - g\circ f$. This Lie algebra is often denoted as $\liegl(V)$, $\liegl_{n}(V)$, or 
 even $\liegl_{n}(k)$ since $V \cong k^n$.
\end{ex}

\begin{ex}
 Given any real Lie group $G$, there is a Lie algebra $\lieg = T_{e}(G)$ corresponding to the tangent bundle on the indentity $e\in G$, where the 
 Lie bracket is given by
 \[
  [X,Y]f = (XY)f - (YX)f
 \]
 for any smooth vector fields $X,Y\in T_{e}(G)$ and any smooth function $f\: G\to \RR$. In the case of a real matrix Lie group $G$, we also can obtain its associated Lie algebra $\lieg$ through the exponential map:
 \[
  \lieg = 
  \left\{
  X \in M^{n\times n}(\RR)\suchthat \forall t \in \RR,\;\; \exp(tX)\in G
  \right\}.
 \]
 In such a case, the Lie bracket is given by the matrix commutator $[X,Y] = XY - YX$. Same construction defines a complex Lie algebra for a complex Lie group $G$.
\end{ex}

\begin{Def}
Let $V$ be an $n$-dimensional vector space over a field $k$. 
 A {\em representation} of a Lie group $G$ is a homomorphism
 $
  \rho\:G\to GL(V) = GL_{n}(k).
 $
 If there exists a proper subspace $W \subset V$ which is invariant for the represeentation $\rho\in \Hom\big(G,GL(V)\big)$, we say that $\rho$ is reducible to $GL(W)$. Otherwise we say that $\rho$ is irreducible. 
\end{Def}

\begin{Def}
A complex Lie group is \emph{reductive} if every representation splits into a direct sum of irreducible representations. 
\end{Def}

Given a compact connected Lie group $K$, its complexification $K_{\mathbb{C}}$ can be proved to be a reductive group.

\subsection[Algebraic varieties]{Algebraic varieties}
\label{ssec:alg-var}

Algebraic varieties are the zero loci of polinomials. Given an algebraically closed field $k$, denote by $\mathbb{A}^{n}_{k}$ the $n$-dimensional affine space over $k$. Let $S$ be a set of polinomials in the ring $k[X_{1}, \ldots, X_{n}]$. An affine algebraic variety is the subset of the affine space where all elements of $S$ take the value zero, this is,
$$V=Z(S)=\{x\in \mathbb{A}^{n}_{k} | f(x)=0,\;\; \forall f\in S\}\; .$$
Sometimes, the definition of affine algebraic variety is restricted to those subsets which are irreducible, meaning that they cannot be described as the union $Z(S_{1})\cup Z(S_{2})$ of two proper subsets defined by the vanishing of polynomials. In these cases, $Z(S)$ is called an algebraic set. For example, $X_{1}^{2}+X_{2}^{2}=1\subset \mathbb{A}^{2}_{\mathbb{C}}$ is the circle contained in the plane, which is an irreducible affine algebraic variety. 

Now define the $n$-dimensional projective space over $k$, $\mathbb{P}_{k}^{n}$. Given a subset of homogeneous polynomials, $S\subset k[X_{0},X_{1}, \ldots, X_{n}]$, a projective variety is a projective algebraic set
$$V=Z(S)=\{x\in \mathbb{P}^{n}_{k} | f(x)=0,\;\; \forall f\in S\}\; ,$$
which is irreducible. As an example, $X_{1}^{2}+X_{2}^{2}=X_{0}^{2}\subset \mathbb{P}^{2}_{\mathbb{C}}$, is the projectivization of the affine circle, which is a projective variety. 
Projetive varieties $X$ can be covered by affine varieties $X\cup U_{i}$, where $U_{i} = \{[X_{0}:\ldots : X_{n}]: X_{i}\neq 0\}$.

Given a projective variety $V\subset \mathbb{P}_{k}^{n}$, we call the \emph{affine cone} of $V$ to the affine algebraic variety $\hat{V}\in k^{n+1}$ resulting of placing a line for each $x\in V$, intersecting in the origin, and yielding a cone such that its transverse slicings recover the projective variety $V$. 

By declaring all algebraic sets $Z(S)$ to be closed, the \emph{Zariski topology} is defined in the projective space and, hence, inherited by any projective variety. A quasi-projective variety is a Zariski open subset of  a projective variety and note that affine varieties are quasi-projective. 

Important examples of algebraic varieties are the {\em Grassmanians}. Given an $n$-dimensional $k$-vector space $V$, the set of all $r$-dimensional vector subspaces of $V$ defines a projective algebraic variety $Gr(r,n)$. Homogeneous coordinates can be given to this variety by the relations of the minors constructed with the generators of subspaces $W\subset V$, which are the Pl\"ucker coordinates embedding $Gr(r,n)$ into a projective space. The Grassmanian $Gr(1,n)$ actually corresponds to the projective space itself.

\subsection[Vector bundles]{Vector bundles}
\label{ssec:v-bles}

Let $X$ be a smooth complex projective variety. Note that, from the analytic point of view, $X$ is also a smooth complex manifold.
Here, we will introduce some basic definitions about vector bundles over $X$. For a more general treatment, the reader may consult \cite{At} or \cite{BT}.

\begin{Def}
 A \emph{holomorphic vector bundle} over $X$ is a smooth manifold $E$ together with a smooth morphism $p\:E\longrightarrow X$ with the following properties:
In a vector bundle, there is an open covering $\{U_{i}\}$ of $X$ such that for every $U_{i}$ there is a biholomorphism $h_{i}$ making the following diagram commutative
$$\xymatrix{
E|_{U_{i}} \ar[r]^{h_{i}} \ar[d]_{p|_{U_{i}}} &
U_{i}\times \mathbb{C}^{r} \ar[ld]^{p_{1}}\\
U_{i}}$$ where $p_{1}$ is projection to the first factor (in other words, $p$ is \emph{locally trivial}). And for every pair $(i,j)$, the composition
$h_{i} \circ h_{j}^{-1}$ is linear on the fibers, i.e.
$$h_{i} \circ h_{j}^{-1}(x,v) = (x,g_{ij}(v))\,,$$ where $g_{ij}:U_i\cap U_{j} \longrightarrow GL(r,\mathbb{C})$ is a holomorphic morphism. The morphisms $g_{ij}$ are called \emph{transition functions}.
For each $x\in X$, the fibers $E_{x} = p^{-1}(x)$ are finite dimensional vector spaces over $\CC$ of dimension $r=\rk(E):= \dim(E_{x})$, called the \emph{rank} of $E$. 
The space $E$ is called the \emph{total space}, the continuous map $p\: E\to X$ is called the \emph{projection map}, and $X$ is called the \emph{base space}.
\end{Def}

Frequently, the total space $E$ denotes the vector bundle altogether, when the base space $X$ and the projection $p$ are clear from the context and no confusion arises.

\begin{Def}
 A \emph{section} of a vector bundle $E$ is a continuous map $s\: X\to E$ such that $(p\circ s)(x) = x$ for all $x\in X$.
\end{Def}

\begin{Def}
 Given two vector bundles $p\: E \to X$ and $q\: F\to X$, a {\em vector bundle homomorphism} is a continuous map $\varphi \: E\to F$ such that $q \circ \varphi = p$, and the restriction to fibers $\varphi \: E_{x} \to F_{x}$ is a linear transformation of vector spaces for each $x\in X$.
The homomorphism $\varphi$ is an {\em isomorphism} if it is bijective and $\varphi^{-1}$ is continuous. We say then that $E$ and $F$ 
 are {\em isomorphic}.
\end{Def}

The usual operations that we carry out on vector spaces and homomorphisms between them, extend naturally to vector bundles, provided that those 
operations are performed fiberwise. For instance, let $p\: E\to X$ and $q\: F\to X$ be two vector bundles over $X$ of ranks
$\rk(E) = n$ and $\rk(F) = m$ respectively. We can define the following bundles:
\begin{enumerate}[i.]
 \item 
 direct sum $E \oplus F$, of rank $n+m$,
 \item 
 tensor product $E \ox F$, of rank $n\cdot m$,
 \item 
 dual bundles $E^{*}$ and $F^{*}$, of ranks $n$ and $m$,
 \item 
 subbundles $F\sse E$ and quotient bundles $E/F$,
 \item 
 exterior powers $\bigwedge^{k}E$ for $k \leq n$ and $\bigwedge^{k}F$ for $k \leq m$,
 \item 
 bundle of homomorphisms $\Hom(E,F)$.
\end{enumerate}
These bundles can be understood by means of operations in the transition functions $g_{ij}$.   For example, transition functions of the direct sum $E\oplus F$ are block diagonal matrices $g^{E}_{ij}\oplus g^{F}_{ij}:U_i\cap U_{j} \longrightarrow GL(n+m,\mathbb{C})$.

\begin{rem}
 Recall that for a subset $F \sse E$ to define a subbundle, the dimension of the fibers $F_{x} = p^{-1}(x)$ must be constant for every
 $x\in X$. Similarly, not every subbundle defines a quotient, $E/F$, it needs to be of constant rank fiberwise as well.
\end{rem}

\begin{ex}
Let $X$ be a smooth complex projective algebraic curve in algebraic geometry, 
(or a compact Riemann surface in complex geometry), of $\dim_{\CC}(X) = 1$.
 The tangent bundle $T(X)$ is an example of a rank $1$ vector bundle where the fibers $E_{x}$ are just copies of $\CC$ encoding all tangent vectors to $X$ at the point $x$. 
 The cotangent bundle $K = K_X = T^{*}(X)$, also called the {\em canonical line bundle} over $X$, is the dual of the tangent bundle. Both are examples of \emph{line bundles}, this is rank $1$ vector bundles. 
\end{ex}

\begin{ex}
Given a projective variety $X$, denote by $\sO_{X}$ its structure sheaf (which is a sheaf of rings, see \cite[Chapter 2]{Ha}) of holomorphic functions on open subsets of $X$. It is the same than a holomorphic line bundle. The line bundle $\sO_{X}(1)$ is the holomorphic bundle whose sections are linear functions in homogeneous coordinates in the projective space. Similarly, $\sO_{X}(d):=\sO_{X}(1)^{\otimes d}$ 
defines a line bundle whose sections are degree $d$ polynomials. 
\end{ex}

Let $X$ be smooth complex projective curve, i.e. a projective variety of complex dimension $1$. Define a \emph{divisor} $D$ in $X$ as a formal finite sum of points 
$$
D=\underset{i}{\sum}n_{i}x_{i}, 
\word{where}
n_{i}\in \mathbb{Z}
\word{and}
x_i \in X.
$$
For each $x\in X$, there exists a meromorphic function $\displaystyle h = \frac{f_{1}\cdot\ldots \cdot f_{l}}{g_{1}\cdot \ldots \cdot g_{s}}$ in a neighborhood $U$ of $x$ whose zeroes are the points of $D\cap U$ with positive coefficient $n_{i}$ and whose poles are those with negative coefficient. Let $\{U_{i}\}$ be a covering of $X$ where, for each $U_{i}$, $D$ has meromorphic equation $h_{i}$. Therefore, in $U_{ij} = U_{i}\cap U_{j}$ the meromorphic function $g_{ij}=\frac {h_{i}}{h_{j}}$ is a unitary element of $\sO^{*}(U_{ij})$ with no zeroes or poles, and it defines transition functions $$g_{ij}:U_{ij} \longrightarrow GL(1,\mathbb{C})$$ of a line bundle. 

We define the \emph{degree} of a divisor $D=\underset{i}{\sum}n_{i}x_{i}$ over a curve $X$ as the sum of the coefficients
$$
\deg(D)
=
\underset{i}{\sum}n_{i}\; .
$$
With the correspondence between divisors and line bundles, we define the degree of a line bundle as the degree of its associated divisor. Roughly speaking, the degree of a line bundle will be the number of zeroes minus the numbers of poles of a rational section (counted with multiplicity).
If $E$ is a vector bundle over $X$ of rank $r$, we define the \emph{determinant line bundle} as the top exterior product $\det(E) = \bigwedge^{r} E$, and define the degree of a vector bundle of rank $r$ as the degree of its determinant line bundle.

The formal definition of the degree of a torsion free sheaf is given in terms of Chern classes. Let $X$ be a smooth complex projective variety of dimension $n$, embedded in a projective space by means of an ample line bundle $\sO_{X}(1)$ corresponding to a divisor $H$.  Given a vector bundle $E$, its Chern classes are denoted by $c_{i}(E)\in H^{2i}(X;\mathbb{Z})$, and define the \emph{degree} of $E$ by 
$$
\deg_{H}(E)=\int_{X} c_{1}(E)\wedge c_{1}\big(\sO_{X}(1)\big)^{n-1}\; .
$$
Note that, if $X$ is an algebraic curve, the degree is the integral of the first Chern class and does not depend on the polarization $\sO_{X}(1)$.

\begin{Def}
Given a vector bundle $E$ over $X$, define its \emph{Euler characteristic} $\chi(E)$ as
$$
\chi(E)
=
\sum_{i\geq 0} (-1)^{i}h^{i}(X;E)\; ,
$$
the alternating sum of the dimensions of the cohomology groups of $E$.
We define the \emph{Hilbert polynomial} of $E$
$$
P_{E}(m)=\chi(E(m))\; ,
$$ 
where $E(m)=E\otimes \sO_{X}(m)$ (called the \emph{twist} of $E$ by $m$) and $\sO_{X}(m)=\sO_{X}(1)^{\otimes m}$.
\end{Def}

\section[Geometric Invariant Theory]{Geometric Invariant Theory}
\label{sec:git}
Let $G$ be a reductive complex Lie group acting on an algebraic variety $X$. The purpose of Geometric Invariant Theory (abbreviated GIT) is to provide a way to define a quotient of $X$ by the action of $G$ with an algebro-geometric structure.  Here we present a sketch of the treatment; for a deeper undestanding and proofs, see \cite{Mu} and the extended version \cite{MFK}.

\subsection{Quotients and the notion of stability}
\label{ssec:git-results}

The problem of taking the quotient of a variety $X$ by the action of a group has been widely studied because it appears everywhere in mathematics. 
In the case when the variety $X$ is affine there is a simpler solution which dates back to Hilbert's $14^{th}$ problem.
Let $A(X)$ denote the coordinate ring of the affine variety $X$. Nagata \cite{Na} proved that if $G$ is reductive,
the ring of invariants $A(X)^{G}$ is finitely generated, hence is the coordinate ring of an affine variety, therefore we can
define the quotient $X/G$ as the affine variety associated to the ring $A(X)^{G}$.

When taking the quotient of a projective variety by a group $G$, there are some issues which have to be taken into account. First,
one has to do with the separatedness of the quotient space and will led us to the definition of $S$-equivalence, or equivalence
of orbits under the action of $G$. Here it is a simple example which shows this feature in the affine case.

\begin{ex}
\label{hyperboles} Consider the action
$$
\xymatrix{ \sigma:\mathbb{C}^{\ast}\times \mathbb{C}^{2}\ar[r] & \mathbb{C}^{2}\\
 (\lambda,(x,y))\ar@{|->}[r] & (\lambda x,\lambda^{-1} y)}
$$
whose orbits are represented in Figure \ref{hypfig}.
\begin{figure}[h]
   \begin{center}
   \includegraphics[width=9cm]{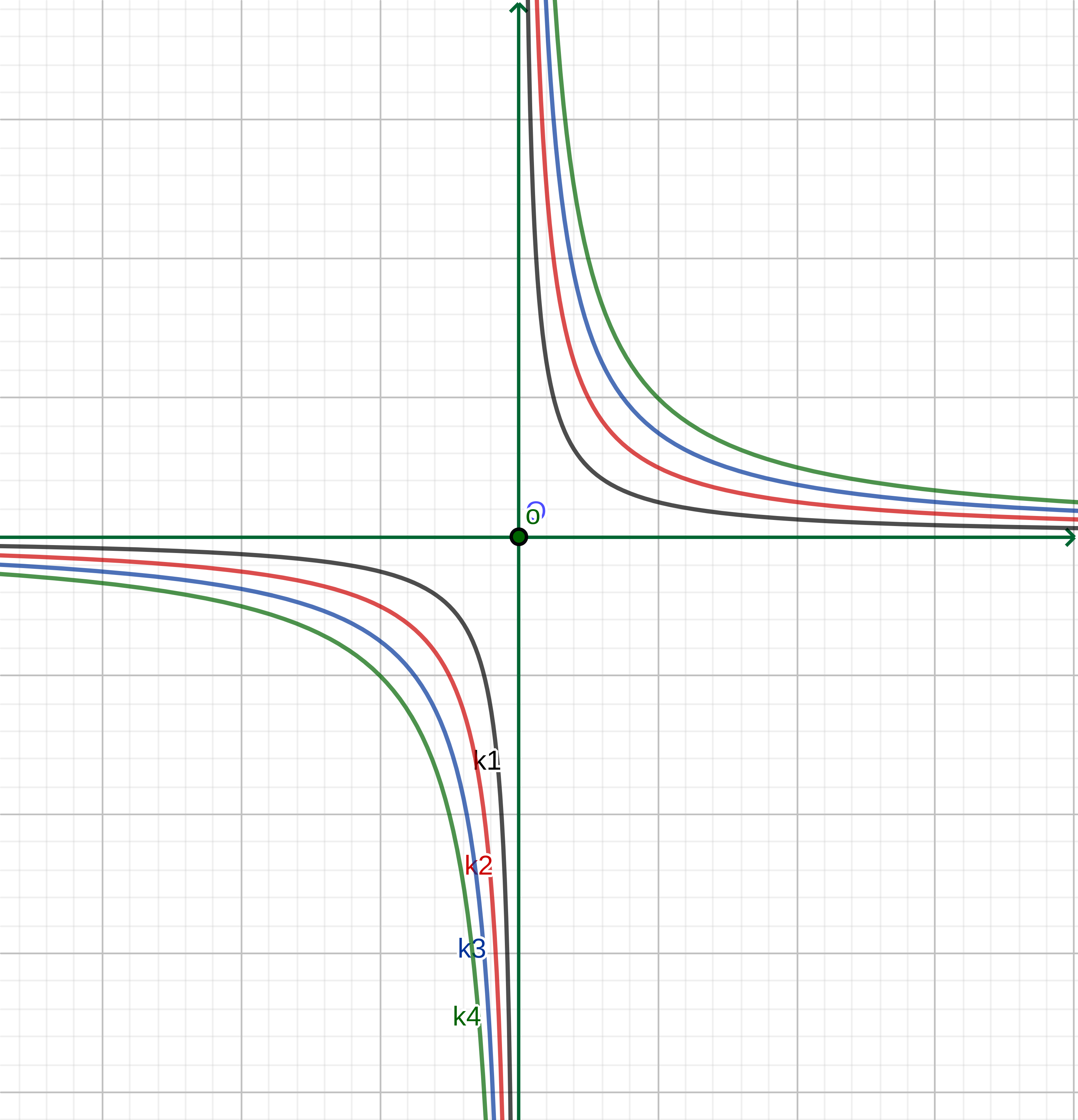}
   \end{center}
\caption{Orbits of the action in Example \ref{hyperboles}.}
\label{hypfig}
\end{figure}
The orbits are the hyperboles $xy=k$, with $k$ a constant, plus three special
orbits, the $x$-axis, the $y$-axis and the origin. Observe that
the origin lies in the closure of the $x$-axis and the $y$-axis.

The coordinate ring of $\mathbb{C}^{2}$ is $\mathbb{C}[X,Y]$, and
the ring of invariants is
$\mathbb{C}[X,Y]^{\mathbb{C}^{\ast}}\simeq \mathbb{C}[XY]\simeq
\mathbb{C}[Z]$. So, the ring of invariants does not distinguish
between the three special orbits, and identifies them in a unique
single point in the quotient space. Hence, the orbit space (the
space where each point corresponds to an orbit) is non
separated, but the quotient space whose ring of functions is
$\mathbb{C}[X,Y]^{\mathbb{C}^{\ast}}\simeq \mathbb{C}[Z]$ is the
affine line, which is separated.
\end{ex}

Once we know how to take quotients of the affine varieties, let us
deal with the projective case. We can guess that, as projective
varieties are given by gluing affine pieces, we can take the
quotient of each affine piece and then glue them together. As we want these
pieces to be respected by the action of $G$ we want them to be
$G$-invariant, hence we are looking for subsets of the form
$$X_{f}=\{x\in X\; |\; f(x)\neq 0\}\; ,$$
which are $G$-invariant or, equivalently, looking for $f\in \mathbb{C}[X_{0},\ldots,X_{n}]$ $G$-invariant.

Now, as the following example shows, note that the action of $G$
on $X$ projective does not determine an action on the graded ring
$\mathbb{C}[X_{0},\ldots,X_{n}]$ (or a quotient of it).

\begin{ex}
Let $\mathbb{C}^{\ast}$ act on $\mathbb{P}^{1}_{\mathbb{C}}$
trivially, i.e. given $\alpha\in \mathbb{C}^{\ast}$, $\alpha\cdot
[x_{0}:x_{1}]= [x_{0}:x_{1}]$. This action is compatible with the
trivial action of $\mathbb{C}^{\ast}$ on $\mathbb{C}[X_{0},X_{1}]$
which acts as $\alpha\cdot f=f$, but it is also compatible with the
action $\alpha\cdot f=\alpha f$ which multiplies each homogeneous polynomial
by the corresponding scalar.
\end{ex}

Hence, we have to linearize the action of $G$ to
$\mathbb{C}^{n+1}$ (called the affine cone of $X$), meaning to give
an action on $\mathbb{C}^{n+1}$ which is the former action of $G$
when restricted to $X$. Once we have this linearization, we can
consider the action on the (graded) coordinate ring of $X$, as we
did in the affine case. We are seeking affine pieces defined as
the complement of the vanishing locus of a $G$-invariant
polynomial, then those points (or orbits) contained on the
vanishing locus of all the $G$-invariant polynomials cannot appear
at any of the affine pieces, hence they cannot be in our quotient.
This motivates the following:

\begin{Def}
A point $x\in X$ is called \emph{GIT semistable} if there exists a $G$-invariant homogeneous polynomial
$f$ of degree $\geq 1$, such that $f(x)\neq 0$. If, moreover, the orbit of
$x$ is closed, it is called \emph{GIT polystable} and if,
furthermore, this closed orbit has the same dimension as $G$ (i.e. if $x$ has finite stabilizer), then $x$ is
called a \emph{GIT stable} point. We say that a closed point of $X$ is \emph{GIT unstable} if it is not GIT
semistable.
\end{Def}

In the previous definition, the idea of semistable points are
those which are separated by homogeneous polynomials, and the
stable ones are those which are infinitesimally separated by
homogeneous polynomials. Indeed, in Example \ref{hyperboles}, all
the orbits $xy=k,\; k\neq 0$, are separated, even infinitesimally,
by the homogeneous polynomial $XY$ (the differential of the
function $XY$ along the transverse direction of the orbits is non
zero) whereas for the three orbits identified (the two axes and
the origin), while they are separated from the other orbits by the
polynomial $XY$, none of them is infinitesimally separated from
the rest. Hence the orbits $xy=k,\; k\neq 0$ are the stable ones
and the other three orbits will define the same point in the quotient,
they will be defined to be equivalent (we will technically say that
they are \emph{S-equivalent}), being the three of them semistable
but not stable and the origin being polystable (the unique closed
orbit in the S-equivalence class).

Note that in this example there are no unstable points, as it will occur in every affine example. Indeed, in affine cases, all points are,
at least, semistable because the constants are always $G$-invariant functions.

\begin{rem}
In general, we consider $X$ embedded in a projective space by the
ample line bundle $\sO_{X}(1)$,
$$X\hookrightarrow \mathbb{P}(H^{0}(\sO_{X}(1))^{\vee})=\mathbb{P}(V)\; .$$
We can see a section $s\in H^{0}(\sO_{X}(m))$ as a
homogeneous polynomial of degree $m$ in $V$. Then, the GIT
unstable points are those for which, for all $m>0$, all
$G$-invariant homogeneous polynomials vanish at that point. This
way, the notion of GIT stability depends on the embedding and the
linearization (i.e. it depends on a line bundle and a lifting of
the action to the total space of this line bundle).
\end{rem}

Mumford \cite{Mu} developed its Geometric Invariant Theory to give a
meaningful geometric structure to the quotient $X/G$. It turns out
that for the semistable orbits we can give a good solution to our
quotient problem. Here we state the technical definition of a good
quotient and the central result of Mumford's GIT.

\begin{Def}
\label{goodquotient} Let $X$ be a projective variety endowed with a
$G$-action. A \emph{good quotient} is a scheme $M$ with a
$G$-invariant morphism $p:X\longrightarrow M$ such that
\begin{enumerate}
\item $p$ is surjective and affine.
\item $p_{*}(\sO^{G}_{X})=\sO_{M}$, where $\sO^{G}_{X}$ is the sheaf of $G$-invariant functions on $X$.
\item If $Z$ is a closed $G$-invariant subset of $X$, then $p(Z)$ is closed in $M$. Furthermore, if $Z_{1}$ and $Z_{2}$ are two closed $G$-invariant subsets of
$X$ with $Z_{1}\cap Z_{2}=\emptyset$, then $p(Z_{1})\cap p(Z_{2})=\emptyset$.
\end{enumerate}
\end{Def}

\begin{thm}\textnormal{\cite[Proposition 1.9, Theorem 1.10]{Mu}}
\label{GIT} Let $X^{ss}$ (respectively, $X^{s}$) be the subset of GIT semistable points (respectively,
GIT stable). Both $X^{ss}$ and $X^{s}$ are open subsets. There is a good quotient $X^{ss}\longrightarrow
X^{ss}/\!\!/G$ (where closed points are in one-to-one correspondence to the orbits of GIT polystable points),
the image $X^{s}/\!\!/G$ of $X^{s}$ is open,
 $X/\!\!/G$ is projective, and the restriction $X^{s}\to X^{s}/\!\!/G$ is a geometric quotient.
\end{thm}

\begin{rem}
The use of double slash $/\!\!/$ in the quotient means that we make
two identifications: one is the identification of the points of
each orbit; the other one is the identification of $S$-equivalent
orbits.
\end{rem}

\begin{rem}
\label{Seq} Two orbits which have non empty intersection will be
called $S$-equivalent and will define the same point in the
quotient. GIT proves that there is only one
closed orbit on each equivalence class (the orbit which is called
polystable). The points of the moduli space are in correspondence
with these distinguished closed orbits, therefore we can say that the moduli that we
obtain classifies polystable points or 
$S$-equivalence classes of points.
\end{rem}

Next we start to analyze the first main example, the construction
of the projective space as a GIT quotient.

\begin{ex}
\label{projective} Let
$\mathbb{C}^{\ast}\; \Circlearrowright\; \mathbb{C}^{n+1}$ be the scalar
action, i.e. $\alpha\cdot (z_{0},\ldots,z_{n})=(\alpha z_{0},\ldots,\alpha z_{n})$,
$\alpha\in \mathbb{C}^{\ast}$, $(z_{0},\ldots,z_{n})\in
\mathbb{C}^{n+1}$. Note that the only invariant functions will be
the constants hence, to have more invariant functions and, then, a
richer quotient space when applying GIT, we will consider invariant
sections of the lifted action by characters (called sometimes in
the literature semi-invariants).

Let $\mathbb{C}^{n+1}\times \mathbb{C}$ be the trivial line bundle
on it and consider different linearizations of the action given by
characters
$$
\chi_{p} \: \xymatrix{\mathbb{C}^{\ast}\ar[r] & \mathbb{C}^{\ast} \\
                \lambda\ar@{|->}[r] & \lambda^{-p}}
$$
such that each linearized action is
$$
\xymatrix{\mathbb{C}^{\ast}\times (\mathbb{C}^{n+1}\times \mathbb{C})\ar[r] & \mathbb{C}^{n+1}\times \mathbb{C} \\
                (\lambda,((z_{0},\ldots,z_{n}),w))\ar@{|->}[r] & ((\lambda z_{0},\ldots,\lambda z_{n}),\lambda^{-p}w)}
$$

If $p>0$ the invariant sections are given by the homogeneous
polynomials
$$f(X_{0},\ldots,X_{n},W)=g_{mp}(X_{0},\ldots,X_{n})\cdot W^{m}\; ,$$
where $g_{mp}$ is a degree $m\cdot p$ homogeneous polynomial on
$(X_{0},\ldots,X_{n})$. The origin will be the unique unstable
point (all $G$-invariant homogeneous polynomials do vanish
simultaneously just at the origin). The semistable locus (indeed
the stable locus, given that all rays are closed orbits in the
semistable locus with maximal dimension) will be
$X^{ss}=\mathbb{C}^{n+1}-\{0\}$. And the quotient, by Theorem
\ref{GIT}, will be a projective variety which we do represent by
$$
X^{ss}/\!\!/\mathbb{C}^{\ast}=\mathbb{P}_{\mathbb{C}}^{n}\; .
$$

If $p<0$, note that there are no invariant sections, hence all orbits are unstable and the quotient is empty.
If $p=0$ the only invariant functions are the constants (it corresponds to the trivial character), hence we cannot separate
any of the orbits from the others and obtain a single point as a quotient.

This example shows how GIT stability of the orbits depends essentially on the different choice of linearization, giving completely different
GIT quotients for different linearizations.

\end{ex}

\subsection{Hilbert-Mumford criterion}
\label{ssec:h-m-criterion}

To determine whether an orbit is GIT stable or unstable we have to
calculate invariant sections or functions. This calculation is
quite involved and dates back to Hilbert. One of Mumford's major
achievements was to give a very simple numerical criterion to
determine GIT stability, called in the literature the
\emph{Hilbert-Mumford criterion}.

It can be proved that a point $x$ is GIT semistable if $0\notin
\overline{G\cdot \hat{x}}$, where $\hat{x}$ lies over $x$ in the
affine cone. Intuitively one direction is clear. Recall that the
GIT unstable points are those for which, for all $m>0$, all
$G$-invariant homogeneous polynomials vanish at that point. As all
homogeneous polynomials (in particular the $G$-invariant ones)
vanish at zero, the points containing zero in the closure of
their orbits will be GIT unstable. The converse can be seen in
\cite[Proposition 4.7]{Ne} or \cite[Proposition 2.2]{Mu}.

The essence of the Hilbert-Mumford criterion is that GIT stability for the whole group $G$ can be checked through $1$-parameter subgroups
$$
\rho:\mathbb{C}^{\ast}\rightarrow G\; ,
$$
stating that we can reach every point in the closure of an orbit
through these $1$-parameter subgroups. Hence, a point is GIT
(semi)stable for the action of $G$ if and only if it is so for the
action of every $1$-parameter subgroup. Then, with the observation
of the previous paragraph, GIT stability measures whether $0$
belongs to the closure of the lifted orbit or not, a belonging which
can be checked through $1$-dimensional paths.

\begin{thm}
\label{HMcrit0} Let $\hat{x}$ be a point in the affine cone over
$X$, lying over $x\in X$.
    \begin{itemize}
    \item $x$ is \emph{semistable} if for all $1$-parameter subgroups $\rho$, $\exists\;\underset{t\rightarrow
    0}{\lim}
    \rho(t)\cdot \hat{x}\neq 0$
        or $\underset{t\rightarrow 0}{\lim} \rho(t)\cdot \hat{x}=\infty$.
    \item $x$ is \emph{polystable} if it is semistable and the orbit of $\hat{x}$ is closed.
    \item $x$ is \emph{stable} if for all $1$-parameter subgroups $\rho$,
    $\underset{t\rightarrow 0}{\lim} \rho(t)\cdot \hat{x}=\infty$ (then the stabilizer of $x$ is finite).
    \item $x$ is \emph{unstable} if there exists a $1$-parameter subgroup $\rho$ such that
    $\underset{t\rightarrow 0}{\lim} \rho(t)\cdot \hat{x}=0$.
    \end{itemize}
\end{thm}

Given $\rho$, a $1$-parameter subgroup of $G$, and given $x\in X$,
we can define $\Phi:\mathbb{C}^{\ast}\longrightarrow X$ by
$\Phi(t)=\rho(t)\cdot x$. We say that $\underset{t\rightarrow 0}{\lim}\;
\rho(t)\cdot x=\infty$, if $\Phi$ cannot be extended to a map
$\widetilde{\Phi}:\mathbb{C}\longrightarrow X$. If $\Phi$ can be
extended, we write $\underset{t\rightarrow 0}{\lim}\; \rho(t)\cdot
x=x_{0}$. The point $x_{0}$ is, clearly, a fixed point of the
action of $\mathbb{C}^{\ast}$ on $X$ induced by $\rho$. Thus,
$\mathbb{C}^{\ast}$ acts on the fiber of the line bundle over
$x_{0}$, say, with weight $\rho_{x}$. One defines the numerical
function
$$\mu(x,\rho):=\rho_{x}\; .$$
We will call this number $\rho_{x}$ the {\em weight} of the
action of $\rho$ over $x$.

The $1$-parameter subgroups induce a linear action of $\mathbb{C}^{\ast}$ in the total space of the line bundle, which
we think as $\mathbb{C}^{n+1}$ for an $n$-dimensional projective variety $X$. By a result of Borel, an action like that can be diagonalized
such that there exists a basis $e_{0},\ldots,e_{n}$ of $\mathbb{C}^{n+1}$ with
$$\rho(t)\cdot \hat{x}=t^{\rho_{i}}\hat{x}_{i}e_{i}\; .$$
Taking into account this, the previous definition of $\mu(\rho,x)$ can be restated as 
$$\mu(x,\rho)=\min\{\rho_{i}\,:\,\hat{x}_{i}\neq 0\}\; .$$

Being defined $\mu(x,\rho)$, we are ready to state the Hilbert-Mumford numerical criterion of GIT stability:

\begin{thm}[{\bf Hilbert-Mumford numerical criterion}]\textnormal{\cite[Theorem~2.1]{Mu}, \cite[Theorem~4.9]{Ne}}
\label{HMcrit} With the previous notations:
    \begin{itemize}
    \item $x$ is \emph{semistable} if for all $1$-parameter subgroups $\rho$, $\mu(x,\rho)\leq 0$.
    \item $x$ is \emph{polystable} if $x$ is semistable and for all $1$-parameter subgroups $\rho$
    such that $\mu(x,\rho) = 0$, $\exists g\in G$ with $x_{0}=g\cdot x$.
    \item $x$ is \emph{stable} if for all $1$-parameter subgroups $\rho$,
    $\mu(x,\rho)<0$.
    \item $x$ is \emph{unstable} if there exists a $1$-parameter subgroup $\rho$ such that
    $\mu(x,\rho)>0$.
    \end{itemize}
\end{thm}

\begin{ex}
 In Example \ref{projective} we can easily check the GIT stability of the orbits by using the numerical Hilbert-Mumford criterion. In
 this case, there is essentially one $1$-parameter subgroup up to rescaling, hence we can directly calculate the minimum
 weight for the action of the group, $\mu(x,\rho)$.

 For all $p$, the action of $\mathbb{C}^{\ast}$ can be extended to the origin in $\mathbb{C}^{n+1}$ which is a fixed point
 for the action. On the fiber over the origin, the action is given by multiplying by $\cdot \lambda^{-p}$, hence for all points
 $x\neq 0$ the minimum weight is $\rho_{x}=-p$ for this ``unique'' $1$-parameter subgroup we are allowed to consider.
 Therefore, by the Hilbert-Mumford criterion, if $p<0$, all points $x\neq 0$ are GIT unstable and, if $p>0$, all points $x\neq 0$ are
 stable. When $p<0$
 it is also clear that the origin is unstable because the minimum weight is $-p$, which is positive. 

However, when $p>0$ we can ``choose another''
 $1$-parameter subgroup (for example with $\cdot \lambda^{2p}$) to obtain a positive weight too, yielding unstability for the point.
 Essentially, the origin is a fixed point for the action and there is no possible linearization making it stable.

 If $p=0$ we have $\rho_{x}=0$ for all points $x$, all orbits are semistable and $S$-equivalent and the only polystable orbit
 is the origin, because it is the limit $x_{0}=0$ not contained in any other orbit but a fixed point.
 \end{ex}

The next example is the fundamental one: the moduli space of binary forms or configurations of $n$ points in the projective line.
It is originally due to Hilbert and it is the starting point for GIT. See \cite{Gi2} for details.

\begin{ex}
\label{npoints}
Let $N$ be an integer and consider the set of all
homogeneous polynomials of degree $N$ in two variables with coefficients in $\mathbb{C}$,
$$
V_{n}=\{f(X,Y)=a_{0}Y^{n}+a_{1}XY^{n-1}+a_{2}X^{2}Y^{n-2}+\cdots+a_{n-1}X^{n-1}Y+a_{n}X^{n}\; |\; a_{i}\in \mathbb{C}\}\; .
$$ 
Let $\mathbb{P}(V_{n})$ be its projectivization. The zeroes of an element $\overline{f}\in \mathbb{P}(V_{n})$ define $n$ points in
$\mathbb{P}_{\mathbb{C}}^{1}$ counted with multiplicity, up to action of the group $G=SL(2,\mathbb{C})$,
$$
\xymatrix{SL(2,\mathbb{C})\times \mathbb{P}(V_{n})\ar[r] & \mathbb{P}(V_{n}) \\
                (g,\overline{f})\ar@{|->}[r] & \overline{f}(g^{-1}(X,Y))}
$$

The orbit space $\mathbb{P}(V_{n})/G$ is not a variety, because it is
not Hausdorff. To see this, let $\overline{f}$ and $\overline{g}$ be represented by
$f=X^{n}$ and $g=X^{n}+X^{n-1}Y$ respectively. The orbits of these two elements are disjoint because the only root of $f$ is
$[0:1]$ counted with multiplicity $n$, and $g$ has two roots,
$[0:1]$ counted with multiplicity $n-1$ and the simple root $[1:-1]$. Let $h_{t}=\left(\begin{array}{cc}
                  t & 0 \\
                    0 & t^{-1} \\ \end{array}
                               \right)$ be a curve of elements in $SL(2,\mathbb{C})$ and define
$$g_{t}:=h_{t}\cdot g=g(h_{t}^{-1}(X,Y))=g(t^{-1}X,tY)=t^{-n}X^{n}+t^{-n+2}X^{n-1}Y\; .$$
For each $t$, $g_{t}$ defines an element in $\mathbb{P}(V_{n})$ which can be represented (by rescalling) by
$\overline{g_{t}}= X^{n}+t^{2}X^{n-1}Y$. Then, note that when $t$ goes to $0$, $\overline{g_{t}}$ tends to
$X^{n}=\overline{f}$, therefore $\overline{f}$ lies in the closure of the orbit of $\overline{g}$ and the orbit space is not
Hausdorff.

In order to construct a GIT quotient we are going to apply the
Hilbert-Mumford criterion in Theorem \ref{HMcrit}. The $1$-parameter subgroups of
$SL(2,\mathbb{C})$ can be diagonalized to be represented by a
diagonal matrix as
$$\rho_{k}(t)=\left(
               \begin{array}{cc}
                 t^{-k} & 0 \\
                 0 & t^{k} \\
               \end{array}
             \right)$$ such that if we write $f(X,Y)=\sum_{i=0}^{n} a_{i}X^{i}Y^{n-i}$ the action of $\rho$ is given by
$$\rho_{k}(t)(f)=f(\rho_{k}(t)^{-1}\cdot (X,Y))=f(t^{k}X,t^{-k}Y)=\sum_{i=0}^{n} a_{i}t^{k(2i-n)}X^{i}Y^{n-i}\; .$$
The limit $f_{0}=\lim_{t\rightarrow 0}\rho_{k}(t)\cdot f$ is equal to the monomial $a_{i_{0}}X^{i_{0}}Y^{n-i_{0}}$, where $i_{0}$ is the
minimum index such that $a_{i}=0$. For example, if $f=XY^{4}+X^{3}Y^{2}+X^{5}$, then $\rho_{k}(t)(f)=t^{-3k}XY^{4}+t^{k}X^{3}Y^{2}+
t^{5k}X^{5}\sim XY^{4}+t^{4k}X^{3}Y^{2}+
t^{8k}X^{5},$ (when considering the projectivization) which tends to $XY^{4}=f_{0}$ when $t$ goes to $0$.

Note that the weight $\rho_{k}$, which acts on the fiber of $f_{0}$, is $\rho_{k,f}=k(2i_{0}-n)$ (in the example, $\rho_{k,f}=-3k$). The Hilbert-Mumford
criterion in Theorem \ref{HMcrit} states that a point $\overline{f}$ is unstable if there exists a $1$-parameter subgroup such that this
weight is positive. Also observe that, up to conjugation in $SL(2,\mathbb{C})$ (or change of homogeneous coordinates $[x:y]$), all $1$-parameter
subgroups are of the form diagonal form $\rho_{k}$ hence classified by the exponent $k$. Therefore, $f$ is unstable if, after a change of coordinates $f(X,Y)=\sum_{i=0}^{n} a_{i}X^{i}Y^{n-i}$, there
exists a $1$-parameter subgroup $\rho_{k}$ such that $k(2i_{0}-n)>0$ where $i_{0}$ is the minimum index such that $a_{i_{0}}\neq 0$.
Given that $k(2i_{0}-n)>0\Leftrightarrow i_{0}>\frac{n}{2}$, it is equivalent to say that $f$ is unstable if and only if $f$ has a root of multiplicity
greater that $n/2$.

In the example of the polynomial $f=XY^{4}+X^{3}Y^{2}+X^{5}$, the weight is
$\rho_{k,f}=-3k<0$, and the lifted orbit
$\rho_{k}(t)(f)=t^{-3k}XY^{4}+t^{k}X^{3}Y^{2}+ t^{5k}X^{5}$ tends
to infinity when $t$ goes to $0$, hence this $1$-parameter
subgroup does not destabilize the point $f$. Indeed, it will occur
the same with all $1$-parameter subgroups as it is easy to check,
because $f$ has no root of multiplicity $\geq 3$. However, the
point $g=X^{3}Y^{2}$ will be acted by $\rho_{k}$ as
$\rho_{k}(t)(g)=t^{k}X^{3}Y^{2}$, which goes to $0$ when $t$ goes
to $0$. Hence, $0$ is in the closure of the lifted orbit and the
weight is $k>0$, then the point is GIT unstable. Indeed $g$ has a
root with multiplicity $3$ (in these coordinates the root is
$[1:0]$).

Observe that, if $n$ is odd, we cannot have $i_{0}=\frac{n}{2}$, hence we cannot have strictly semistable points and all the GIT semistable
points will be GIT stable.

If $n$ is even, we can observe the $S$-equivalence phenomenon. Let
$n=4$ and consider the points $f=X^{2}Y^{2}+X^{3}Y+X^{4}$ and
$g=X^{2}Y^{2}$. By the same argument that we used to show that the
orbit space is not Hausdorff, it is clear that $f$ (with $2$ roots
equal and the other two different) and $g$ (with two roots
pairwise equal) do not lie in the same orbit but $g$ lies in the closure
of the orbit of $f$. Hence the $2$ points are $S$-equivalent. To
determine which one is the only polystable orbit within this
equivalence class we can use a $1$-parameter subgroup of type
$\rho_{2}$ which acts on the fiber of the limit point (common to
$f$ and $g$ and, indeed, equal to $g$) with weight zero, and conclude
that $g$ is the polystable orbit.
\end{ex}

\begin{rem}
The moduli space of configurations of $n$ points in the projective line is the same that the moduli space of $n$-gons (c.f. \cite{KM}) if
we consider the isomorphism $\mathbb{P}_{\mathbb{C}}^{1}\simeq S^{2}$ and see points in $\mathbb{P}_{\mathbb{C}}^{1}$ as
length unit vectors. A configuration of points will be unstable if there is a point with multiplicity more than half the points, the
same way a polygon will be unstable if there is any of the vectors repeated more that half times. It can be shown that
an unstable polygon does not close, in the sense that, after any change of coordinates by $SL(2,\mathbb{C})$, the sum of
the vectors is not zero.
\end{rem}

\subsection{Symplectic stability}
\label{ssec:sympl-stability}

In this subsection we will sketch the symplectic reduction procedure, giving another perspective of the stability picture. The Kempf-Ness
theorem will be the link in between the two sides. Let us begin by reviewing the basics about symplectic geometry. 

Let $(X,\omega)$ be a \emph{symplectic manifold}, where $X$ is a smooth manifold and $\omega\in \Omega^{2}(X)$ is a closed non-degenerate
two form (called a \emph{symplectic form}). Two symplectic varieties $(X_{1},\omega_{1})$ and $(X_{2},\omega_{2})$ are 
\emph{symplectomorphic} if there exists a diffeomorphism $\varphi:X_{1}\rightarrow X_{2}$ such that $\varphi^{\ast}\omega_{2}=\omega_{1}$. 
By Darboux's theorem every symplectic manifold is locally symplectomorphic to $\mathbb{R}^{2n}$ equipped with the standard symplectic
$2$-form $\sum_{i=1}^{n}dq_{i}\wedge dp_{i}$. 

Given a symplectic manifold $(X,\omega)$, let $\Symp (X,\omega)\subset \Diff(X)$ be the group of symplectomorphisms and let 
$\Vect^{s}(X)\subset \Vect(X)$ be the Lie subalgebra of symplectic vector fields $v\in \Vect(X)$ such that $\mathcal{L}_{v}\omega=
d(\iota_{v}\omega)0$. 
Given a smooth function $H\in\mathcal{C}^{\infty}(X,\mathbb{R})$, it defines a symplectic vector field $\xi_{H}$ by $\iota_{\xi_{H}}\omega=dH$. 
Observe that the image of $\mathcal{C}^{\infty}(X,\mathbb{R})$ lies in the subalgebra $\Vect^{s}(X)$ of symplectic vector fields. 
In local Darboux coordinates, $\xi_{H}$ is given by 
$$\xi_{H}=\sum_{i=1}^{n}\frac{\partial H}{\partial p_{i}}\frac{\partial}{\partial q_{i}}-\frac{\partial H}{\partial q_{i}}
\frac{\partial}{\partial p_{i}}\; ,$$
from which we can see, by remembering the Hamilton equations, how symplectic geometry gives the natural framework 
for mechanics. We call $d\mathcal{C}^{\infty}(X,\mathbb{R})=\Vect^{H}(X)$ the hamiltonian vector fields. Given that the kernel of $d$ are the 
constant functions, the Lie algebra of the hamiltonian automorphisms is $\mathcal{C}^{\infty}(X,\mathbb{R})/\mathbb{R}$.

Let $K$ be a compact connected Lie group acting on a symplectic manifold $(X,\omega)$. We say that the action is symplectic if it 
preserves the symplectic form, i.e. $k_{X}\in \Symp(X,\omega),\; \forall k\in K$. We say that the action is \emph{hamiltonian}
if the map $\mathfrak{k}\rightarrow \Vect(X)$ (which sends an element $\xi\in\mathfrak{k}=Lie(K)$ to the corresponding vector 
field in $X$) lifts, equivariantly by the action of $K$, to a hamiltonian vector field $\xi_{H}$, $H\in \mathcal{C}^{\infty}(X,\mathbb{R})$ such that $\iota_{\xi_{H}}\omega=dH$.
In this case, we can define a moment map 
$$\mu:X\rightarrow \mathfrak{k}^{\ast}$$
by the condition $\iota_{\xi}\omega=-d\langle\mu,\xi\rangle$, $\forall \xi\in\mathfrak{k}$. Given that the Lie algebra of the hamiltonian 
automorphisms is $\mathcal{C}^{\infty}(X,\mathbb{R})/\mathbb{R}$, we can choose each element $\xi_{H}$ up to a constant; hence the lifting
condition means that we choose these constants in such a way that $\mu$ is $K$-equivariant (by the coadjoint action on the 
right hand side). Therefore, given a hamiltonian $K$-action, the moment map is unique up to the addition of a central element
of $\mathfrak{k}^{\ast}$. 

In the following, let $X\subset \mathbb{P}_{\mathbb{C}}^{n}$ be a projective variety with an action of a compact connected
Lie group $K$, whose complexified group is $G$ (which is, hence, reductive). For simplicity, consider that $G\subset GL(n+1,\mathbb{C})$
and $K\subset U(n+1)$. Suppose that $K$ acts on $\mathbb{P}_{\mathbb{C}}^{n}$ by preserving the almost-complex structure $J$
and the Fubiny-Study metric $g$, hence $K$ preserves the natural symplectic structure $\omega=g(\cdot,J\cdot)$. In this case
there is a natural moment map which, for $K=U(n)$ and identifying the Lie algebra $\mathfrak{u}(n)$ with its dual
(via the inner product $\langle A,B \rangle=\trace(A^{\ast}B)$), is given by
\begin{equation}
 \label{momentmap}
\mu:\mathbb{P}_{\mathbb{C}}^{n}\rightarrow \mathfrak{u}(n)^{\ast}\;,\;\mu(z)=\frac{i}{2}zz^{\ast}\; ,
\end{equation}
up to addition of a central element which in this case is a constant. When we have a diagonal action on a product of
symplectic varieties it can be proved that the moment map is the sum of the respective moment maps. 

\begin{rem}
The different moment maps for a given action correspond with the different polarizarions and linearizations of the action from 
the Geometric Invariant Theory side. If the symplectic form $\omega$ is integral, meaning that its cohomology 
class lies in $H^{2}(X,\mathbb{Z})/torsion\leq H^{2}(X,\mathbb{R})$, then $2\pi i\omega$ is the curvature of an hermitian line bundle
$L$ with a unitary connection, and the isometries of $L$ preserving the connection cover the hamiltonian authomorphisms on $X$. 

In the projective case, the cohomology class is integral, hence we can develop this \emph{prequantization} to restrict to a discrete
number of different moment maps, associated to the GIT linearizations. 
\end{rem}

In the symplectic setting we state the following notion of stability.

\begin{Def}
\label{mustab} Let $(X,\omega)$ be a projective variety with the symplectic form coming from the Fubini-Studi metric, endowed with a hamiltonian $K$-action. 
Let $\mu$ be a moment map for 
this action. Let $x$ be a point of $X$ and let us denote by $G\cdot x$ its orbit by the complexified group $G=K^{\mathbb{C}}$.
    \begin{itemize}
    \item $x$ is \emph{$\mu$-semistable} if $\overline{G\cdot x}\cap \mu^{-1}(0)\neq \emptyset$.
    \item $x$ is \emph{$\mu$-polystable} if $G\cdot x\cap \mu^{-1}(0)\neq \emptyset$.
    \item $x$ is \emph{$\mu$-stable} if $x$ is $\mu$-polystable and, in addition, the stabilizer of $x$ under $G$ is finite.
    \item $x$ is \emph{$\mu$-unstable} if $\overline{G\cdot x}\cap \mu^{-1}(0)=\emptyset$.
    \end{itemize}
\end{Def}

The notions of GIT stability and $\mu$-stability will be equivalent by the Kempf-Ness theorem.

\begin{thm}[\textbf{Kempf-Ness Theorem }\cite{KN}]
Let $(X,\omega)$ be a projective variety with the symplectic form coming from the Fubini-Studi metric, endowed with a hamiltonian $K$-action.
Let $\mu$ be a moment map for this action. A $G$-orbit is GIT polystable if and only if it contains a zero of the moment map.
A $G$-orbit is GIT semistable if and only if its closure contains a zero of the moment map, and this zero lies in the unique GIT
polystable orbit in the closure of the original orbit. 
\end{thm}

We will make some considerations to sketch the proof of the Kempf-Ness theorem. 

Let $(X,L=\sO_{X}(1))$ be a projective polarized variety and choose an hermitian metric on $L$ inducing
a connection with curvature $2\pi i\omega$. Lift a point $x\in X$ to $\hat{x}\in L_{x}^{-1}$ and consider the functional norm
$\|\hat{x}\|$. If $X\subset \mathbb{P}(H^{0}(L)^{\ast})$ and we consider a metric in $H^{0}(L)^{\ast}$, it induces a metric in the total space
of $L^{-1}$ where $\|\hat{x}\|$ is the norm in the vector space where the affine cone $\hat{X}$ lives. 

For each $\hat{x}$, define the Kempf-Ness function
\begin{equation}
\psi_{\hat{x}}:\mathfrak{k}\rightarrow \mathbb{R},\;\;\; \xi\mapsto \frac{\log\|\exp(i\xi)\hat{x}\|^{2}}{2}\; .
\end{equation}
The $1$-parameter subgroups encoding GIT stability by the Hilbert-Mumford criterion can be thought as 
different directions in the $G$-orbit, hence different elements of the Lie algebra $Lie(G)=\mathfrak{g}=
\mathfrak{k}\oplus i\mathfrak{k}$. To study how this function varies along $1$-parameter subgroups we calculate
$$\partial_{\lambda}\psi_{\hat{x}}(\xi)=\frac{d}{dt}|_{t=0}\frac{\log\|\exp(i(\xi+t\lambda))\hat{x}\|^{2}}{2}=$$
$$\frac{\langle i\lambda\exp(i\xi)\tilde{x},\exp(i\xi)\hat{x}\rangle}
{\langle \exp(i\xi)\hat{x},\exp(i\xi)\hat{x}\rangle}=2\mu((\exp i\xi)x)(\lambda)\; ,$$
which can be expressed by saying that the Kempf-Ness function is an integral of the moment map. If we calculate the second derivative, we obtain
$$\partial_{\nu}\partial_{\lambda}\psi_{\tilde{x}}(\xi)=2\langle \mathcal{L}_{J\nu}\mu((\exp i\xi)x),\lambda\rangle=$$
$$(\omega(\lambda,J\nu))(\exp(i\xi)x)=g(\lambda,\nu)(\exp(i\xi)x)\; ,$$
which is non negative, since $g$ is a Riemannian metric.

Hence, the Kempf-Ness function is convex, attaining a minimum at the zeroes of the function $\mu((\exp i\xi)x)$ which are
the zeroes of the moment map. This way, $x$ is $\mu$-polystable if and only if $\psi_{\hat{x}}$ attains a minimum. If
 the Kempf-Ness function is bounded from below it does not necessarily attain a 
minimum but, if it does asymptotically, it means that the closure of the $G$-orbit of the point contains a zero of the moment map and the 
point is $\mu$-semistable. 

The unstable points will be those for which the Kempf-Ness function is not bounded from below or, equivalently, 
the orbit under the complexified group does not intersect the zeroes of the moment map. The GIT unstable points are those $x$ for which $0\in \overline{G\cdot \hat{x}}$, where $\hat{x}$ lies over $x$ in the affine cone. From the definition 
of the Kempf-Ness function $\psi_{\hat{x}}$ in terms of the logarithm, $0\in \overline{G\cdot \hat{x}}$ will be equivalent to $\psi_{\hat{x}}$ not to be 
bounded by below, which is equivalent to the $\mu$-unstability of $x$. 

From this, the symplectic quotient construction is due to Meyer and Marsden-Weinstein \cite{Me, MW}:

\begin{thm}
Let $(X,\omega)$ be a symplectic manifold endowed with a hamiltonian action of a compact connected Lie group $K$. If $\mu$ is a moment
map for this action, and $K$ acts freely and properly on $\mu^{-1}(0)$, the quotient $(\mu^{-1}(0)/K,\omega_{0})$ is a smooth
symplectic manifold with $i^{\ast}\omega=p^{\ast}\omega_{0}$, where $i:\mu^{-1}(0)\rightarrow X$ is the inclusion and
$p:\mu^{-1}(0)\rightarrow \mu^{-1}(0)/K$ the projection, respectively. 
\end{thm}

By the Kempf-Ness theorem, we will have the following bijection relating the GIT and the symplectic quotients, which is indeed 
an isomorphism:
$$
\mu^{-1}(0)/K\simeq X^{ps}/G=X^{ss}/\!\!/G\; .
$$


\subsection{Examples}
\label{ssec:examples-stability}
Next, we will calculate the moment map for the examples studied from the algebraic setting and check that the Kempf-Ness theorem 
holds in these cases.

\begin{ex}
\label{projectivesym}
Let us go back to Example \ref{projective}. The compact group in this case is $K=U(1)\subset GL(1,\mathbb{C})=\mathbb{C}^{\ast}$. 
In this case the different moment maps are given by (c.f. (\ref{momentmap}) and \cite{Wo})
$$\xymatrix{\mu: \mathbb{C}^{n+1}\ar[r] & \mathfrak{k}^{\ast}=\mathfrak{u}(1)^{\ast}\simeq \mathbb{R}\\
(z_{0},\ldots,z_{n})\ar@{|->}[r]  & \frac{i}{2}(|z_{0}|^{2}+\cdots|z_{n}|^{2}-a)}$$
where $a$ comes from a central element of $\mathfrak{k}^{\ast}=\mathfrak{u}^{\ast}$, which in this case is 
any real number. If we 
add the condition that the lifted action of $\mathfrak{u}$ descends to an action of the group $K=U(1)$ on the trivial line bundle
we have that $a\in \mathbb{Z}$. The different $a\in \mathbb{Z}$ correspond to the integers $p$ of the different characters 
in Example \ref{projective}. 

If $a<0$ there are no $\mathbb{C}^{\ast}$-orbits intersecting $\mu^{-1}(0)$, not even in the closure, hence all points are
$\mu$-unstable as well as they were GIT unstable. 

If $a=0$, the origin in $\mathbb{C}^{n+1}$ is $\mu$-polystable because its orbit intersects $\mu^{-1}(0)$ and all the other orbits
are $\mu$-semistable but not $\mu$-polystable because their closures intersect $\mu^{-1}(0)$. The origin is in the closure of all orbits, hence it
is the unique polystable point in the unique $S$-equivalence class. Therefore, the symplectic quotient is again a single point. 

If $a>0$, the origin is $\mu$-unstable because its orbit does not intersect $\mu^{-1}(0)$. All the other orbits intersect $\mu^{-1}(0)$
at some $(z_{0},\ldots,z_{n})$ such that $\sum_{i=0}^{n}|z_{i}|^{2}=a$, hence all rays are $\mu$-polystable (indeed $\mu$-stable)
and the quotient
is the expected projective space $\mathbb{P}_{\mathbb{C}}^{n}$. 
\end{ex}

\begin{ex}
 \label{npointssym}
Now we recall the classification of configurations of $n$ points in $\mathbb{P}_{\mathbb{C}}^{1}$, from Example \ref{npoints}.
Identify each $f\in \mathbb{P}(V_{n})$ with the set of its $n$ zeroes counted with multiplicity and, by the isomorphism
$\mathbb{P}_{\mathbb{C}}^{1}\simeq S^{2}$, identify them with $n$ vectors in the unit sphere. 
The compact group now is $SO(3,\mathbb{R})\subset SL(2,\mathbb{C})$, acting diagonally on $(S^{2})^{n}$ by rotations.
The Lie algebra of $SO(3,\mathbb{R})$ is $\mathfrak{so}(3,\mathbb{R})\simeq \mathbb{R}^{3}$ and the moment map in this case is just the sum of the inclusions
of each vector in $\mathbb{R}^{3}$, hence given by (c.f. \cite{Wo})
$$\xymatrix{\mu: (S^{2})^{n}\ar[r] & \mathfrak{su}(2)^{\ast}\simeq \mathbb{R}^{3}\\
(v_{1},\ldots,v_{n})\ar@{|->}[r]  & v_{1}+\cdots v_{n}}$$
Then, a configuration of points will be $\mu$-semistable if and only if the associated $n$-tuple of vectors 
$(v_{1},\ldots,v_{n})$ (up to action of
the complexified group $SL(2,\mathbb{C})$), verify $\sum_{i=1}^{n}v_{i}=0$, which is the equivalent to say that a ``polygon closes'', 
identifying this problem with the moduli space of polygons (see \cite{KM}). 

Since the Kempf-Ness theorem asserts that $\mu$-stability is equal to GIT stability, this means that a configuration of $n$ points
in $\mathbb{P}_{\mathbb{C}}^{1}$ can be moved, by an element of $SL(2,\mathbb{C})$, such that
the corresponding $n$-tuple of vectors in $S^{2}$ (counted with multiplicity) have center of mass the origin, if and only if
there is no point with multiplicity greater than half the total, which means that the point is semistable.

In the case $n$ is even, we can have a point with multiplicity exactly half the total (recall that this meant the point is
GIT semistable but not stable). The polynomials $f=X^{2}Y^{2}+XY^{3}+Y^{4}$ and $g=X^{2}Y^{2}$ verify that
$g\in \overline{SL(2,\mathbb{C})\cdot f}$. The polynomial $g$ defines a configuration with only two points,
each of them with the same multiplicity equal to half the total. For example the points $[1:0]$ and $[0:1]$, which in $S^{2}$ can be thought as the vectors 
$(0,0,1)$ and $(0,0,-1)$. Then, $g$ is the only polystable orbit in the closure of the orbit of $f$ which defines a (degenerate) configuration
of vectors in the unit sphere with center of mass the origin, therefore $\mu^{-1}(0)\cap G\cdot g\neq\emptyset$ and 
$\mu^{-1}(0)\cap \overline{G\cdot f} \neq \emptyset$ but $\mu^{-1}(0)\cap G\cdot f=\emptyset$, meaning that $f$ is $\mu$-semistable
but not $\mu$-polystable and $g$ is $\mu$-polystable. 
By visualizing polygons, this situation in general corresponds to the degenerate polygon with 
$n/2$ vectors equal to $v$ and the other $n/2$ equal to $-v$ 
lying on a line, which can only appear for $n$ even. This limit point corresponds to 
the polystable orbit with stabilizer $\mathbb{C}^{\ast}$. 
\end{ex}

\begin{ex}
\label{grass}
We will obtain the Grassmannian as a GIT quotient and as a symplectic quotient. 

Let $SL(r,\mathbb{C})\Circlearrowright \Hom(\mathbb{C}^{r},\mathbb{C}^{n})$, $r<n$, be the group action such that $A\cdot g^{-1}$ for
$A\in \Hom(\mathbb{C}^{r},\mathbb{C}^{n})$, $g\in SL(r,\mathbb{C})$, and linearize the induced action on the projectivized 
vector space $\mathbb{P}(\Hom(\mathbb{C}^{r},\mathbb{C}^{n}))$ to the tautological line bundle by 
$$
g\cdot ([A],\lambda)=\lambda A\cdot g^{-1}\; ,
$$ 
where $\lambda$ is an element of the fiber of the tautological line bundle lying over $[A]$.
The points of the Grassmannian of $r$-planes in $\mathbb{C}^{n}$ will correspond to injective homomorphisms from $\mathbb{C}^{r}$ to 
$\mathbb{C}^{n}$, up to change of basis. This change of basis is encoded by considering the projectivized $\mathbb{P}(\Hom(\mathbb{C}^{r},\mathbb{C}^{n}))$
(two linear maps differing by multiplication of a scalar define the same $r$-plane) and by the action of $SL(r,\mathbb{C})$ (changes
of frame with determinant $1$). 
Hence, let us prove that $[A]\in \mathbb{P}(\Hom(\mathbb{C}^{r},\mathbb{C}^{n}))$ is GIT stable if and only if
$A\in \Hom(\mathbb{C}^{r},\mathbb{C}^{n})$ has rank $r$. 

If $\rk A< r$, pick a basis $\{v_{1},\ldots,v_{r}\}$ of $\mathbb{C}^{r}$ such that $v_{1}\in \Ker A$. Choose a $1$-parameter subgroup
$\rho$ adapted to the basis such that it has the diagonal form 
$$\left(
               \begin{array}{cccc}
                 t^{r-1} & 0 &  \cdots & 0\\
                 0 & t^{-1} & \cdots & 0\\
                \vdots  &    & \ddots & \vdots \\
                 0 & 0 & \cdots &  t^{-1} 
               \end{array}
             \right)$$
Then, $$A\cdot \rho^{-1}=\left(\begin{array}{cccc}
                 0 & \ast& \cdots & \ast\\
                 0 & \ast & \cdots & \ast\\
                \vdots  & \vdots  &  & \vdots\\
                 0 & \ast & \cdots & \ast 
               \end{array} \right)\cdot
               \left(\begin{array}{cccc}
                 t^{1-r} & 0 &  \cdots & 0\\
                 0 & t^{1} & \cdots & 0\\
                \vdots  &   & \ddots & \vdots\\
                 0 & 0 & \cdots &  t^{1} 
               \end{array}
             \right)=\left(\begin{array}{cccc}
                 0 & t\cdot\ast & \cdots & t\cdot \ast\\
                 0 & t\cdot\ast & \cdots & t\cdot\ast\\
                \vdots  & \vdots  &  & \vdots\\
                 0 & t\cdot\ast & \cdots & t\cdot\ast 
               \end{array} \right)=t\cdot A\; ,$$
hence $\rho$ fixes $[A]\in\mathbb{P}(\Hom(\mathbb{C}^{r},\mathbb{C}^{n}))$ and acts on the fiber $\mathbb{C}\cdot A$ as $\cdot t$, this
is with weight $\rho_{A}=1>0$, therefore $[A]$ is GIT unstable. 

Conversely, if $A$ has full rank, up to action of $SL(2,\mathbb{C})$, there exists a splitting $\mathbb{C}^{n}\simeq 
\mathbb{C}^{r}\oplus \mathbb{C}^{n-r}$ where $A$ is the inclusion of the first factor in this splitting. Given $\rho$, a $1$-parameter
subgroup of $SL(2,\mathbb{C})$, we assume that we can choose a basis which both diagonalizes $\rho$ and agrees with the splitting.
Then, $\rho$ is 
$$\left(         \begin{array}{ccccc}
                 t^{\lambda_{1}} & 0 & 0 & \cdots & 0\\
                 0 & t^{\lambda_{2}} & 0 & \cdots & 0\\
                 0 & 0 & t^{\lambda_{3}} & \cdots & 0\\
                 \vdots  & \vdots  &   & \ddots & \vdots\\
                 0 & 0 & 0 & \cdots & t^{\lambda_{r}} 
               \end{array}
             \right)$$
             and assume further that $\lambda_{1}\geq \lambda_{2}\geq \cdots \geq \lambda_{r}$, with $\sum_{i=1}^{r}\lambda_{i}=0$. Note that,
             by rescalling in $\mathbb{P}(\Hom(\mathbb{C}^{r},\mathbb{C}^{n}))$, the action of $\rho$ in $[A]$, i.e. 
             $[A]\cdot \rho^{-1}$, is the same that $[A]\cdot \rho^{-1}\cdot t^{\lambda_{1}}$. Then, the diagonal of
             $\rho^{-1}\cdot \lambda_{1}$ is $(1,\ldots,1,t^{-\lambda_{i}+\lambda_{1}},\ldots,t^{-\lambda_{r}+\lambda_{1}})$, where all
             $-\lambda_{i}+\lambda_{1}>0$ (if $\rho$ is not trivial). When we take the limit $t\rightarrow 0$, $A$ tends to $A_{0}$
             where $A_{0}$ represents
             the inclusion of $\mathbb{C}^{p}$ as the first $p$ vectors of the basis in $\mathbb{C}^{n}$ ($p$ is the number
             of $1$'s in the diagonal of $\rho^{-1}\cdot \lambda_{1}$, equal to the number of exponents $\lambda_{1}$ in $\rho$). 
             Finally, the weight of $\rho$ in the fiber over the limit point $A_{0}$ is $-\lambda_{1}<0$, and $[A]$ is GIT stable. 
             
Equivalently, from the symplectic point of view, we have the action of the unitary group $U(r)\subset GL(r,\mathbb{C})$ acting on 
$\Hom(\mathbb{C}^{r},\mathbb{C}^{n})$ the same way. By considering the inner product $\langle A,B \rangle=\trace(A^{\ast}B)$ which identifies 
$\mathfrak{u}^{\ast}(n)$ with $\mathfrak{u}(n)$, a moment map for the action is (c.f. (\ref{momentmap}))
$$\xymatrix{\mu: \Hom(\mathbb{C}^{r},\mathbb{C}^{n})\ar[r] & \mathfrak{u}(r)^{\ast}\\
A\ar@{|->}[r]  & \frac{i}{2}(A^{\ast}A-Id)}$$
Hence, $\mu^{-1}(0)$ are those matrices such that, up to action of $GL(r,\mathbb{C})$, verify $A^{\ast}A=Id$, which is to say 
that a linear map is congruent by $GL(r,\mathbb{C})$ to an isometric embedding if and only if it is injective. 

In general, we could have added a central element (in this case a scalar) to the moment map to get $\mu(A)=\frac{i}{2}(A^{\ast}A-a\cdot Id)$. If $a>0$ we 
obtain the same result. If $a=0$ the quotient is a single point and if $a<0$ all points are $\mu$-unstable. This corresponds
to the different linearizations in the GIT problem. 
\end{ex}

\subsection{Maximal unstability}
\label{ssec:maximal-unstability}

After studying the relation between GIT stability and symplectic stabilily by the Kempf-Ness theorem, in this section we will
focus on the unstable locus. We will classify the unstable points by degrees of unstability and will check that this notion 
agrees when considered from both points of view. 

The moment map $\mu:X\rightarrow \mathfrak{k}^{\ast}$ 
is invariant by the adjoint action of the compact
group $K$ but not by the action of its complexified group $G=K^{\mathbb{C}}$. If we choose an inner product 
$\langle \cdot,\cdot \rangle$ 
in $\mathfrak{k}$, invariant by $K$, we can identify $\mathfrak{k}^{\ast}$ with $\mathfrak{k}$ and 
define the function $\|\mu\|:X\rightarrow \mathbb{R}$ by $\|\mu(x)\|=\langle \mu(x),\mu(x)\rangle$, to which 
we will refer as the \emph{moment map square}. Recall 
that the Kempf-Ness function is an integral of the moment map. The $\mu$-unstable points are those $x$ for which $\mu(g\cdot x)$  does not
achieve zero as a limit point, for $g\in G$, hence the Kempf-Ness function for these points is unbounded. 

Define the function $\Omega_{x}(g)=\|\mu(g\cdot x)\|$, $g\in G$. The function $\Omega_{x}$ is a Morse-Bott function
and it takes some infimum value $m_{x}\geq 0$ at the critical set. The idea is that there exists 
a direction of maximal descense for the 
negative gradient flow of the Kempf-Ness function, directions thought as cosets in $G/K$, minimizing the moment map square, i.e. 
the function $\Omega_{x}$ (c.f. \cite{Ki} and \cite{GRS}). Then, the $G$-orbit of a $\mu$-unstable point $x$ does not achieve 
$\Omega_{x}^{-1}(0)$ but it achieves, in their closure, $\Omega_{x}^{-1}(m_{x})$
for some positive number $m$ (c.f. Moment limit theorem \cite[Theorem 6.4]{GRS} 
and Generalized Kempf Existence Theorem \cite[Theorem 11.1]{GRS}). Of course, for the $\mu$-semistable ones this 
infimum $m_{x}$ is zero.

From the algebraic point of view, recall that a point $x$ is GIT unstable if there exists a $1$-parameter subgroup $\rho$ such that
the weight $\rho_{x}$ is positive (recall that the number $\rho_{x}$ is the weight that $\rho$ is acting with on 
the fiber of the fixed limit point of $\rho(t)$ when $t$ goes to zero). Having chosen the inner product $\langle \cdot,\cdot \rangle$ 
in $\mathfrak{k}$, it extends uniquely to an inner product in $G$. Considering the $1$-parameter subgroups as directions given
by elements in the Lie algebra $\mathfrak{g}=\Lie (G)$, it makes sense to define the norm $\|\rho\|$ of a $1$-parameter subgroup
and define the function $\Phi_{x}(\rho)=\rho_{x}/\|\rho\|$. If $x$ is GIT unstable, there exists $\rho$ such that
$\Phi_{x}(\rho)>0$. The result in \cite{Ke} asserts that the supremum of the function $\Phi_{x}$ is attained at some unique $\rho$
(up to conjugation by the parabolic subgroup of $G$ defined by $\rho$), hence
there exists a unique $1$-parameter subgroup maximizing the Hilbert-Mumford criterion, or giving the maximal way to destabilize
a GIT unstable point. The norm in the denominator serves to calibrate this maximal degree of unstability when rescalling 
(i.e. multiplying the exponents of the $1$-parameter subgroups by a scalar).  

The principal result in \cite{GRS} (c.f. \cite[Theorem 13.1]{GRS}) shows that, for $x$ an unstable point,
$$\sup_{\rho\in \mathfrak{g}}\Phi_{x}(\rho)=\sup_{\rho\in \mathfrak{g}}\frac{\rho_{x}}{\|\rho\|}=m_{x}=\inf_{g\in G}\Omega_{x}=
\inf_{g\in G}\|\mu(g\cdot x)\|\; ,$$
this is, the weight of the $1$-parameter subgroup which maximally destabilizes a GIT unstable point $x$ (after
normalization) is the infimum of the moment map square over the $G$-orbit of a $\mu$-unstable point.

\begin{ex}
Let us go back to Example \ref{npoints}, the configurations of points in $\mathbb{P}^{1}_{\mathbb{C}}$. The group $SL(2,\mathbb{C})$ is simple, 
then there is only one invariant inner product up to multiplying by a scalar, say the Killing norm. Then, we can choose 
$\langle \cdot,\cdot \rangle$ such that the associated norm verifies $$\left\|\left(\begin{array}{cc}
                 t^{-k} & 0\\
                 0 & t^{k}  
               \end{array}\right)\right\|=k\; .$$
We did calculate in Example \ref{npoints} that the weight of a $1$-parameter subgroup $\rho_{k}$ which has exponents $-k$ and $k$ in 
its diagonal form is $\rho_{f}=k(2i_{0}-n)$ where, recall that $i_{0}$ is the maximum number of points in $\mathbb{P}^{1}_{\mathbb{C}}$
which are equal. It is clear that 
$$\sup_{\rho\in \mathfrak{g}}\Phi_{f}(\rho)=\sup_{\rho\in \mathfrak{g}}\frac{\mu(f,\rho)}{\|\rho\|}=\frac{k(2i_{0}-n)}{k}=2i_{0}-n\; ,$$
which is a positive number if $f$ is unstable. 

Now, from the symplectic point of view, recall that we associate to each point in $\mathbb{P}^{1}_{\mathbb{C}}$ a vector in 
$S^{2}$ and the moment map is given by $\mu(x)=v_{1}+\cdots +v_{n}\in \mathbb{R}^{3}$, after identifying 
$\mathfrak{so}(3,\mathbb{R})^{\ast}\simeq \mathbb{R}^{3}$. The norm chosen in $\mathfrak{so}(3,\mathbb{R})$ can be identified with the usual norm in 
$\mathbb{R}^{3}$. 

Suppose that $x$ is an unstable configuration, hence it defines $i_{0}>\frac{n}{2}$ identical vectors in $S^{2}$. By changing the coordinates
in $\mathbb{P}^{1}_{\mathbb{C}}$, we can consider that the configuration is given by a binary form 
$$f=a_{n-i_{0}}X^{n-i_{0}}Y^{i_{0}}+a_{n-i_{0}+1}X^{n-i_{0}-1}Y^{i_{0}+1}+\cdots +a_{n-1}XY^{n-1}+a_{n}Y^{n}\; ,$$
which can be moved in its $G$-orbit by elements $g_{t}=\left(\begin{array}{cc}
                 t & 0\\
                 0 & t^{-1}  
               \end{array}\right)$ to obtain
$$g_{t}\cdot f=fg_{t}^{-1}=f(t^{-1}X,tY)=$$
$$t^{2i_{0}-n}a_{n-i_{0}}X^{n-i_{0}}Y^{i_{0}}+
t^{2i_{0}-n+2}a_{n-i_{0}+1}X^{n-i_{0}-1}Y^{i_{0}+1}+\cdots +t^{n-2}a_{n-1}XY^{n-1}+t^{n}a_{n}Y^{n}\; .$$
We can multiply it by $t^{-2i_{0}+n-2}$ and still define the same form $\overline{f}$ in the projective space,
$$a_{n-i_{0}}X^{n-i_{0}}Y^{i_{0}}+
t^{2}a_{n-i_{0}+1}X^{n-i_{0}-1}Y^{i_{0}+1}+\cdots +t^{2n-2i_{0}-2}a_{n-1}XY^{n-1}+t^{2n-2i_{0}}a_{n}Y^{n}\; ,$$ 
which tends to $f_{0}=a_{n-i_{0}}X^{n-i_{0}}Y^{i_{0}}$ when $t$ goes to $0$. The zeroes of $\overline{f_{0}}$ are $[1:0]$ 
with multiplicity $i_{0}$ and $[0:1]$ with multiplicity $n-i_{0}$ and, when considering a isomorphism $\mathbb{P}^{1}_{\mathbb{C}}\simeq S^{2}$, we
can associate the roots to the vectors $(0,0,1)$ and $(0,0,-1)$ in $S^{2}$. Hence, the calculation of the infimum of the moment map square is
$$\inf_{g\in G}\Omega_{f}=
\inf_{g\in G}\|\mu(g\cdot f)\|\leq \inf_{t}\|\mu(g_{t}\cdot f)\|=$$
$$|\sum_{i_{0}}(0,0,1)+\sum_{n-i_{0}}(0,0,-1)|=|\sum_{2i_{0}-n}(0,0,1)|=2i_{0}-n=:m_{f}\; ,$$
and it is clear that the value obtained is indeed the infimum, because the best we can do in order
to get the infimum, once we have $i_{0}$ identical vectors in 
$S^{2}$, is to dispose the rest (up to action of $SL(2,\mathbb{C})$) in the opposite direction, which we did by the curve of 
elements $g_{t}\in G$. 

As we observe, 
$$\sup_{\rho\in \mathfrak{g}}\Phi_{f}(\rho)=2i_{0}-n=\inf_{g\in G}\Omega_{f}\; ,$$
therefore there are different levels of unstability, indexed by the numbers $m_{f}=2i_{0}-n$, corresponding to binary forms with
different number of identical roots, or to vectors in $S^{2}$ which do not close to form a polygon because they have different numbers of identical vectors, in all cases more than half of them.            
\end{ex}

\begin{ex}
 Now we recall Example \ref{grass}. 

 Let $A\in \Hom(\mathbb{C}^{r},\mathbb{C}^{n})$ of rank $m<r$, hence $[A]$ is GIT unstable. Following the argument in the example, there exists a basis
 $\{v_{1},\ldots,v_{r}\}$ of $\mathbb{C}^{r}$ such that $L\{v_{1},\ldots,v_{r-m}\}=\Ker A$. The different $1$-parameter subgroups
$\rho$, adapted to the basis in such a way they take the diagonal form, are given by 
$$\left(
               \begin{array}{cccc}
                 t^{\lambda_{1}} & 0 &  \cdots & 0\\
                 0 & t^{\lambda_{2}} & \cdots & 0\\
                \vdots  &    & \ddots & \vdots \\
                 0 & 0 & \cdots &  t^{\lambda_{r}} 
               \end{array}
             \right)$$
where we impose the convention $\lambda_{1}\geq\lambda_{2}\geq\cdots\geq\lambda_{r}$. Then, 
$$A\cdot \rho^{-1}=\left(\begin{array}{ccccc}
                 0 & \cdots & \ast & \cdots & \ast\\
                 \vdots & \ddots & \vdots &  & \vdots\\
                0  & \cdots  & \ast  & \cdots & \ast\\
                \vdots & & \vdots & & \vdots\\
                 0 & \cdots & \ast & \cdots & \ast\\
                 
               \end{array} \right)\cdot
               \left(\begin{array}{cccc}
                 t^{-\lambda_{1}} & 0 &  \cdots & 0\\
                 0 & t^{-\lambda_{2}} & \cdots & 0\\
                \vdots  &    & \ddots & \vdots \\
                 0 & 0 & \cdots &  t^{-\lambda_{r}} 
               \end{array}
             \right)=$$
             
             $$\left(\begin{array}{ccccc}
                 0 & \cdots & t^{-\lambda_{r-m+1}}\cdot \ast & \cdots & t^{-\lambda_{r}}\cdot\ast\\
                 \vdots & \vdots & \vdots &  & \vdots\\
                0  & \cdots  & t^{-\lambda_{r-m+1}}\cdot \ast  & \cdots & t^{-\lambda_{r}}\cdot\ast\\\
                \vdots & & \vdots & & \vdots\\
                 0 & \cdots & t^{-\lambda_{r-m+1}}\cdot \ast & \cdots & t^{-\lambda_{r}}\cdot\ast\
               \end{array} \right)$$
Hence, we observe that the weight $\mu([A],\rho)$ of the Hilbert-Mumford criterion, i.e. the minimal exponent multiplying a non-zero coordinate,
is $-\lambda_{r-m+1}$. Therefore, in order to maximize this weight, keeping the condition that $\rho\in SL(r,\mathbb{C})$ hence all exponents sum up to $0$, 
the maximal $1$-parameter subgroups will be of the form
$$\left(
               \begin{array}{cccccc}
                 t^{m} & 0 &  & \cdots & & 0\\
                 0 & \ddots & & & & 0\\
                  &  & t^{m}  &  & & \\
                \vdots & & & t^{m-r} & & \vdots\\
                  & & &  & \ddots & \\
                0 & 0 &  & \cdots & & t^{m-r} 
               \end{array}
             \right)$$
 where the exponent $m$ is repeated $r-m$ times and the exponent $m-r$ is repeated $m$ times. Then it is clear that for these $1$-parameter subgroups
we have $\mu([A],\rho)=r-m$. Note that we could have achieved the same maximal result by multiplying the exponents $m$ and $m-r$ by the same positive 
constant hence, up to rescalling, the maximal weight will remain $r-m$. In other words,
$$\sup_{\rho\in \mathfrak{g}}\Phi_{[A]}(\rho)=\sup_{\rho\in \mathfrak{g}}\frac{\mu([A],\rho)}{\|\rho\|}=r-m\; .$$

From the symplectic side, recall that the moment map was given by $\mu(A)=\frac{i}{2}(A^{\ast}A-Id)$. Having chosen the invariant product in $\mathfrak{u}(n)$
given by $\trace(A^{\ast}B)$, the moment map square is given by
$$\|\mu(A)\|=\trace((A^{\ast}A-Id)^{\ast}(A^{\ast}A-Id))=\trace((A^{\ast}A-Id)^{2})$$
up to a constant (related with the rescalling of the norm discussed before from the GIT point of view). By an element of $SL(r,\mathbb{C})$ (or
by change of basis) we can suppose that $A^{\ast}A$ is a matrix with a diagonal block which is the idendity (of size the rank of
$A$) and zeroes elsewhere. Therefore, it is clear that 
$$\inf_{g\in G}\Omega_{A}=\inf_{g\in G}\|\mu(g\cdot A)\|=$$
$$\trace\left(\left(\begin{array}{cccccc}
                 0 & 0 &  & \cdots & & 0\\
                 0 & \ddots & & & & 0\\
                  &  & 0  &  & & \\
                \vdots & & & 1 & & \vdots\\
                  & & &  & \ddots & \\
                0 & 0 &  & \cdots & & 1 
               \end{array}\right)-Id\right)^{2}=\trace\left(\begin{array}{cccccc}
                 1 & 0 &  & \cdots & & 0\\
                 0 & \ddots & & & & 0\\
                  &  & 1  &  & & \\
                \vdots & & & 0 & & \vdots\\
                  & & &  & \ddots & \\
                0 & 0 &  & \cdots & & 0 
               \end{array}\right)=r-m\; ,$$
which is equal to the quantity $\sup_{\rho\in \mathfrak{g}}\Phi_{[A]}(\rho)$. Hence, the different unstability levels are indexed by the complementary of
the rank of $A$, being $m=r$ the case where the supremum and the infimum, respectively, achieve zero, as it has to be in the stable case. 
\end{ex}

\section{Moduli Space of vector bundles}
\label{sec:ms}

The problem of classifying vector bundles is a very central story in geometry since the 1960s, with strong connections with other areas of mathematics and physics. The statement is to find a geometric structure with good properties (such as an algebraic variety), where each point corresponds to a holomorphic structure in a given smooth vector bundle. 

When trying to perform this, the objects to classify happen to have different groups of automorphisms. This turns out to be a main issue because the moduli space attempts to collect all structures, identifying the ones which are mathematically identical, this is, modding out by their automorphisms. Therefore, this prevents us from solving the problem of finding a moduli space for all vector bundles. 

However, following the idea of Grothendieck, it is often the case that one can add a piece of data to the objects to classify, such that the only automorphism of these objects endowed with the additional data is the identity. The piece of data we add is encoded as the action of certain group in a parameter space. Geometric Invariant Theory, which was developed in \cite{Mu} for this particular purpose as its main application, gives the solution to the quotient by the action of that group on the objects, removing the additional data, and yielding a moduli space as a GIT quotient.

The notion of stability that we define for vector bundles will distinguish between stable bundles, those with the smallest possible automorphism group and for which the solution of the moduli problem is the best possible (a fine moduli space parametrizing isomorphism classes), semistable (but non-stable) bundles having a coarse moduli space (parametrizing S-equivalence classes) and unstable bundles left out of the classical moduli problem and that need to be dealt with by means of the Harder-Narasimhan filtration.  

\subsection{GIT construction of the moduli space}
\label{ssec:GITms}

Let $X$ be a smooth algebraic projective curve. This is the same as a Riemann surface (a compact topological surface with a holomorphic structure), and let $g$ be its genus. Suppose $X$ embedded in a projective space with a very ample line bundle $H = \sO_{X}(1)$, called a polarization of $X$.

Let $E$ be a holomorphic vector bundle of rank $r$, degree $d$ and fixed determinant line bundle $\det E\simeq L$. Define the \emph{slope} of a vector bundle by $\mu(E)=\frac{d}{r}$.

\begin{Def}
A vector bundle $E$ is said to be \emph{semistable} if for every proper subbundle $0\neq E' \subsetneq E$ we have $$\mu(E') \leq \mu(E)\; . $$
A vector bundle is said to be \emph{stable} if the inequality is strict for every proper subbundle. A vector bundle which is non-semistable will be called \emph{unstable}. A vector bundle $E$ is \emph{polystable} if it is isomorphic to a direct sum of stable bundles 
$$E \simeq E_{1}\oplus E_{2} \oplus \cdots \oplus E_{l}$$
of the same slope $\mu(E_{i})$.
\end{Def}

Note that the degree and, therefore, the notion of stability, does not depend on the polarization of the curve $X$. This is not true for higher dimensional varieties $X$. 

Let m be an integer. A vector bundle $E$ over $X$ is \emph{$m$-regular} if
$H^{i}(E(m-i)) = 0$, $\forall i>0$. If $E$ is $m$-regular, then $E(m):=E\otimes \sO_{X}(m)$ is generated by global sections (the evaluation map is surjective). By the Vanishing Theorem of Serre, for each vector bundle $E$, there exists an integer $m$ (depending on $E$) such that $E$ is $m$-regular. 

It can be proved (see \cite{NS,Se} and \cite{Ma} in higher dimension) that all semistable vector bundles of rank $r$ and degree $d$ are bounded, meaning that they are parametrized by a finite-type scheme whose ranks and degrees of their subbundles and quotients are bounded. Then, it can be shown that, for a bounded family, we can choose an integer $m$ uniformly in Serre's theorem. 

With the choice of $m$, the dimension of the space of global sections of the twisted bundle $E(m)$ is
$$h^{0}(E(m))=\chi(E(m))=\deg(E(m))+r(1-g)=d+rm+r(1-g)=:N\; ,$$
a linear polynomial in $m$. 
Let $V$ be an $N$-dimensional complex vector space. Choose an isomorphism 
$\alpha : V\simeq H^{0}(E(m))$. By composing with the evaluation map we obtain a surjection
$$V\otimes \sO_{X} \simeq H^{0}(E(m))\otimes \sO_{X}  \twoheadrightarrow E(m)\; .$$
Taking cohomology, we get a homomorphism between vector spaces
$$U = H^{0}(U\otimes \sO_{X})\rightarrow H^{0}(E(m))$$ and, taking the $r$-exterior power, 
$$\bigwedge^{r}V\rightarrow \bigwedge^{r} H^{0}(E(m))\rightarrow H^{0}(\bigwedge^{r}(E(m)))= $$
$$
H^{0}(\bigwedge^{r}E\otimes \sO_{X}(rm))\simeq H^{0}(L\otimes \sO_{X}(rm)) =: W
$$
where the isomorphism comes from $\beta: \bigwedge^{r} E =\det E \rightarrow L$, and two of these isomorphisms $\beta$ differ by a scalar. Therefore, a point $(E,\alpha)$  provides an element in $\Hom(\bigwedge^{r} V, W)$ and, because of the choice of $\beta$, a well defined element in $\mathbb{P}(\Hom(\bigwedge^{r} V, W))$.
Then we have all semistable bundles, together with a choice of isomorphism $\alpha$, inside the so-called Quot-scheme (or scheme of quotients) of Grothendieck
$$Z := \{(E,\alpha)\} \subset Quot_{V,m,X,r,d}\hookrightarrow \mathbb{P}(\Hom(\bigwedge^{r}V,W))\; ,$$
parametrizing quotients $V\otimes \sO_{X}(-m) \twoheadrightarrow E$. To obtain a moduli space we will take the GIT quotient of this subvariety $Z$ by the action encodiing the changes of isomorphism $\alpha$, this is the action of $GL(V)=GL(N,\mathbb{C})$ although, because of the projectivity, it is enough to take the quotient by $SL(V)$. Then, the moduli space of semistable vector bundles will be the GIT quotient $\mathcal{M}_{r,L} = Z/\!\!/SL(V)$.

By GIT results (c.f. Theorem \ref{GIT}), there exists a good quotient of the GIT semistable points $Z^{ss}/\!\!/SL(V)$. Hilbert-Mumford criterion (c.f. Theorem \ref{HMcrit}) states that GIT stability can be checked by $1$-parameter subgroups $\rho:\mathbb{C}^{*}\longrightarrow SL(N,\mathbb{C})$. Once fixed m, denote $V=H^{0}(E(m))$ and let $V^{i}$ be each one of the eigenspaces of the diagonalization of $\rho$, where $V=\bigoplus_{i}V^{i}$ and $V_{i}=\bigoplus_{j=1}^{i}V^{j}$, $V^{i}=V_{i}/V_{i-1}$. Let $E_{i}(m)$ be the sheaf generated by the sections of $V_{i}$, with rank $r_{i}:=\rk E_{i}(m)$, and denote by $E^{i}(m)=E_{i}(m)/E_{i-1}(m)$, with rank $r^{i}:=\rk E^{i}(m)$. 

Let us compute the weight of the action in the limit point of the $1$-parameter subgroup, $\mu(\rho, x)$, where $x$ is the point in the Quot scheme, therefore in the projective space, corresponding to $E$. We obtain
$$\mu(\rho,E)=\dim V\cdot \sum_{n\in\mathbb{Z}}n\cdot r^{n}=\sum_{n\in\mathbb{Z}}n(r^{n}\cdot \dim V-r\cdot \dim V^{n})=-\sum_{n\in\mathbb{Z}}(r_{n}\cdot \dim V-r\cdot \dim V_{n})\; ,$$ 
and, written in terms of the exponents $\rho_{i}$ of the diagonalization of $\rho$, get
$$\mu(\rho,x)=\sum_{i=1}^{N}\rho_{i}(-\dim V^{i}\cdot r+\dim V\cdot r^{i})\; .$$
Observe how the final sum is taken over a finite number of non-zero terms.
From here we can see that, $E$ is semistable (resp. stable), if and only if,
$$\dim (V')\cdot r\underset{(<)}{\leq} \dim V\cdot r'\; ,$$ for all $V'\subsetneq V$, $E'(m)$ the subsheaf generated by $V'$ and $r'=\rk(E'(m))$.

Recall that $\dim V=h^{0}(E(m))=\chi(E(m))$ (because of the $m$-regularity of $E$) and note that $\dim V'<h^{0}(E'(m))$, the generated bundle having more sections in general. Then, if $E$ were a GIT unstable point, there would exists $V'\subset V$ such that
$$\dfrac{h^{0}(E'(m))}{r'}>\dfrac{\dim V'}{r'}>\dfrac{\dim V}{r}=\dfrac{h^{0}(E(m))}{r}\; .$$ Since we have, by Riemann-Roch theorem,
$$\dfrac{h^{0}(E(m))}{r}=\dfrac{d+rm+r(1-g)}{r}=\dfrac{d}{r}+(1-g)+m=\mu(E)+(1-g)+m\; ,$$
and similar for $E'(m)$, the inequality turns out to be 
$$\mu(E')>\mu(E)\; ,$$
recovering the definition of non-stability for vector bundles. This shows that an unstable vector bundle yields a GIT unstable point and viceversa, semistable bundles correspond to GIT semistable points. 

Therefore $Z^{ss}$ matches, indeed, the semistable bundles that we want to classify, and the GIT quotient 
$$\mathcal{M}_{r,L} = Z^{ss}/\!\!/SL(V)$$
is the moduli space of semistable vector bundles of rank $r$ and determinant $L$.
As GIT says, this moduli space is a good quotient where points correspond to S-equivalence classes of semistable vector bundles, each class containing a unique closed orbit, being the polystable representative. The stable bundles $Z^{s}/SL(V)$ yield a good quotient which is a quasi projective variety.

\subsection{Harder-Narasimhan filtration}
\label{ssec:HN}

Unstable bundles fall outside the solution to the moduli problem. However, to each of them, there is attached a canonical filtration called the \emph{Harder-Narasimhan filtration} (see \cite{HN} and \cite[Section 1.3]{HL}) which exhibits unstable bundles as extensions of semistable ones.  

\begin{thm}
\label{HNdef}
Let $E$ be a vector bundle of rank $r$. There exists a unique filtration, called the \emph{Harder-Narasimhan filtration} for $E$,
$$\left\{ 0\right\}=E_{0}\subsetneq E_{1}\subsetneq E_{2}\subsetneq ...\subsetneq E_{t-1}\subsetneq E_{t}=E$$
satisfying
\begin{itemize}
\item The slopes of the quotients $E^{i}:=E_{i}/E_{i-1}$ verify
$$\mu(E^{1})>\mu(E^{2})>\mu(E^{3})>...>\mu(E^{t-1})>\mu(E^{t})\; .$$
\item The quotients $E^{i}, \forall i\in \left\{1,...,t\right\}$ are semistable.
\end{itemize}
\end{thm}

As we discussed in the construction of the moduli space, unstable bundles have the biggest possible automorphism group, coming from the different semistable blocks of their Harder-Narasimhan filtration and the homomorphisms between them. On the other hand, stable bundles are much easier because their are \emph{simple} objects whose automorphisms are just the scalars. In between these two, semistable but non-stable bundles can be studied by means of their Jordan-H\"older fitration, whose graded object captures this automorphism group. This will lead us to the notion of $S$-equivalent, the right property to identify semistable points in the moduli space. 

\begin{thm}
\label{JN}
Let $E$ be a semistable vector bundle of rank $r$. There exists a (non unique in general) filtration, called the \emph{Jordan-H\"older filtration} for $E$,
$$\left\{ 0\right\}=E_{0}\subsetneq E_{1}\subsetneq E_{2}\subsetneq ...\subsetneq E_{l-1}\subsetneq E_{l}=E$$
satisfying
\begin{itemize}
\item The slopes of the quotients $E^{i}:=E_{i}/E_{i-1}$ are equal
$$
\mu(E^{1})=\mu(E^{2})=\mu(E^{3})=...=\mu(E^{l-1})=\mu(E^{l})\; .
$$
\item The quotients $E^{i}, \forall i\in \left\{1,...,l\right\}$ are stable.
\end{itemize}
The filtration is unique in the sense that the graded objects $gr(E)\simeq \bigoplus_{i=1}^{l} E^{i}$
of two different Jordan-H\"older filtrations of $E$ are isomorphic. 
\end{thm}

\begin{Def}
Let $E$, $E'$, such that the graded objects of any of their Jordan-H\"older filtrations are isomorphic, $gr(E)\simeq gr(E')$. Then $E$ and $E'$ are said to be \emph{$S$-equivalent}.
\end{Def}

If $E$ and $E'$ are $S$-equivalent bundles, their orbits in the GIT construction of the moduli space are S-equivalent in the sense of GIT, and they correspond to the same point in the moduli space. Vector bundles $E$ which are already isomorphic to the graded object of their Jordan-H\"older filtrations, this is
$$E \simeq E_{1}\oplus E_{2} \oplus \cdots \oplus E_{l}$$
where all slopes $\mu(E_{i})$ are equal and the factors are stable bundles, are precisely the polystable vector bundles  corresponding to GIT polystable orbits in the GIT construction, the only closed orbits in each $S$-equivalence class. This completes the classification of points in the moduli problem.

\subsection{Analytical construction of the moduli space of vector bundles}

After Narasimhan and Seshadri \cite{NS} algebraic construction of the moduli space of holomorphic vector bundles, relating them to the unitary representations, and the GIT of Mumford \cite{Mu, MFK}, providing a compactification of the moduli problem, Atiyah and Bott \cite{AB}, build an equivalent classification by using Morse theory. 

Let $X$ be a smooth Riemann surface (equivalent to a smooth complex projective curve) and let $\sE\to X$ a complex smooth bundle. This is, a vector bundle whose fibers are complex vector spaces, but its transitions functions are complex differentiable not necessarily holomorphic. 

Denote by $\sA^{0,1}(\sE)$ the space of holomorphic structures on $\sE$. Given that the rank $r$ and the degree $d$ of all these structures are the same (once defined $\mathcal{E}$, these are topological invariants), we can denote it by $\sA^{0,1}(r,d)$. It happens that the quotient $\sA^{0,1}(r,d) / \gG^{\CC}$
  by the gauge group\footnote{The gauge group of a principal $G$-bundle $P$ is defined as the space of sections of the adjoint bundle $\Ad P = P\times_{G} G$. It accounts for the idea of an automorphism group of $P$, invariant by the action fo the group $G$. For details, see \cite{AB}.} $\gG^{\CC}$ is not Hausdorff so, in order to have an interesting quotient, we need to impose the stability condition. 

 \begin{Def} 
   The \emph{moduli space of (polystable) vector bundles} $\sE \to X$ is defined as the quotient
   \[
   \sN(r,d) = \sA_{ps}^{0,1}(r,d) / \gG^{\CC} 
   \]
where the complex \emph{gauge group} $\gG^{\CC}$ acts on the subspace
$$
 \sA_{ps}^{0,1}(r,d) = \{\sA^{0,1}(\sE)\: \sE \text{ is polystable}\}\; .$$

Respectively, \emph{the moduli space of stable vector bundles} is defined as the quotient
   \[
   \sN_{s}(r,d) = \sA_{s}^{0,1}(r,d) / \gG^{\CC}
   \subseteq
   \sN(r,d),
   \]
   where
 $
 \sA_{s}^{0,1}(r,d) = \{\sA^{0,1}(\sE)\: \sE \text{ is stable}\}.
 $
 \end{Def}
 
This moduli space $\sN_{s}(r,d)$ is equivalent to the one defined in Section \ref{ssec:GITms}, $\mathcal{M}_{r,L}$, where the degree of $L$ is precisely $d$.


In the coprime case, where there are no strictly semistable objects, already Narasimhan and Seshadri compute the dimension of the moduli space. 
 
 \begin{thm}\textnormal{\cite{NS}}
 \label{thm:dimvectorbundles}
  If $\GCD(r,d) = 1$, then $\sA_{ps}^{0,1} = \sA_{s}^{0,1}$
  and $\sN(r,d) = \sN_{s}(r,d)$ is a smooth complex projective variety of dimension
  \[
  \dim_{\CC}\big(\sN(r,d)\big) = r^2 (g - 1) + 1.
  \]
 \end{thm}

\section[Moduli space of Higgs bundles]{Moduli space of Higgs bundles} 
\label{sec:Higgs}

Higgs bundles were introduced by Hitchin in $1987$ to study the Yang-Mills equations over Riemann surfaces. A Higgs bundle over a Riemann surface $X$ is a pair $(E,\varphi)$ where $E$ is a holomorphic bundle over $X$ and $\varphi$ is a section of the bundle $E\otimes K$, $K$ being the canonical bundle of $X$. They represent the mathematical formulation of the scalar field that, in the \emph{Higgs mechanism}, interacts with gauge bosons making them behave as if they carry mass. This theory intends to explain the symmetry break in particle physics. 

Similarly to holomorphic bundles, by introducing a stability condition, a moduli space of polystable HIggs bundles can be constructed, yielding a complex algebraic variety. Original work \cite{Hi1} for rank $2$ bundles, shows the Hitchin-Kobayashi correspondence where polystable Higgs bundles correspond with certain reduction of the structure group $G=\SL(2,\mathbb{C})$ (in the language of metrics and connections), satisfying generalized Yang-Mills equations. 
This proof is generalized by Simpson \cite{Si1, Si2} and Garc\'ia-Prada, Gothen and Mundet i Riera \cite{G-PGM} for other groups. In this way, the moduli space of Higgs bundles turns out to be homeomorphic to the Dolbeault moduli space, classifying  operators $\overline{\partial}$ in a smooth bundle, equivalent by gauge transformations. 

The other fundamental piece in the Higgs bundle universe is its relationship with the representations of the fundamental group of $X$, $\rho:\pi_{1}(X)\rightarrow G$. Works of Simpson, Corlette and Donaldson, among others, what is known as the non-abelian Hodge theory, provides a homeomorphism of the moduli space of Higgs bundles with the moduli of representations $\rho$ which are reductive, and extends the the theorem of Narasimhan-Seshadri \cite{NS}, adding the Higgs field on the bundle side, and the trip from the real compact form ($U(n)\subset GL(n,\mathbb{C})$ originally) to the group $G$. In this way, it is completed a correspondence between three moduli spaces coming algebra, geometry, topology and physics.

\subsection[Hitchin construction]{Hitchin construction}
\label{ssec:Hitchin-construction}

Let $X$ be a Riemann surface and let $\mathcal{E}$ a complex differentiable (or smooth) bundle over $X$. 
Hitchin (c.f. \cite{Hi1}) establishes a reduction of the \emph{Yang-Mills self-duality equations} (SDE) from $\RR^4$ to $\RR^2 \simeq \CC$:

 \begin{equation} 
 \left\{
 \begin{array}{c c c}
  F_A + [\varphi, \varphi^{*}] & = & 0 \\
                         &   &   \\
  \bar{\partial}_A \varphi  & = & 0,
 \end{array}
 \right.
 \label{eq:YM} 
 \end{equation}
where $\varphi \in \Omega^{1,0}\big(X, \End(\sE)\big)$ is the \emph{Higgs field} and $F_A$ is the curvature of a connection $d_A$ which is compatible with the holomorphic structure of the holomorphic bundle $E = (\sE, \bar{\partial}_A)$, and $\sE$ has rank $\rk(\sE)=2$ and degree $\deg(\sE)=1$. Here, $\varphi^{*}$ denotes the adjoint of $\varphi$ with respect to the hermitian metric on $E$, and $[\cdot,\cdot]$ denotes the natural extension of the Lie bracket to Lie algebra-valued forms. For more details, the reader may see Hitchin \cite{Hi1}.

 The set of pairs which are solutions of SDE
 $$
 \beta(\sE) := 
 \{
 (\bar{\partial}_A, \varphi)| \text{ solution of }(\ref{eq:YM})
 \}
 \subseteq
 \sA^{0,1}(\sE) \times \Omega^{1,0}(X,\End(\sE))
 $$
and the collection
 $$
 \beta_{ps}(2,1) := 
 \{
 \beta(\sE)| \sE \text{ polystable, } \rk(\sE) = 2, \deg(\sE) = 1 
 \},
 $$
allows Hitchin to construct the Moduli space of polystable solutions to SDE~(\ref{eq:YM})
  $$
  \sM^{YM}(2,1) = \beta_{ps}(2,1) / \gG^{\CC},
  $$
and
  $$
   \sM^{YM}_{s}(2,1) = \beta_{s}(2,1) / \gG^{\CC}
   \subseteq
   \sM^{YM}(2,1),
  $$
the moduli space of stable solutions.

\begin{rem}
 Since $\GCD(2,1) = 1$, then $\sA_{ps}^{0,1} = \sA_{s}^{0,1}$, this is, there are no strictly semistable holomorphic structures and the moduli spaces coincide
  $$
  \sM^{YM}(2,1) = \sM^{YM}_{s}(2,1).
  $$
\end{rem}

The work of HItchin \cite{Hi1, Hi2} also presents an alternative algebro-geometric construction. It resembles on the extension of the stability notion to Higgs bundles. Recall that a \emph{Higgs bundle} over $X$ is a pair $(E, \varphi)$ where 
 $E \to X$ is a holomorphic vector bundle and the Higgs field is $\varphi\in H^{0}(E\otimes K)$.
The stability condition is analogous to the one for vector bundles, 
but checking  just on subbundles preserved by $\varphi$.

\begin{Def}
 A subbundle $F \sse E$ is said to be $\varphi$-{\em invariant} 
 if $\varphi(F) \sse F \ox K$. A Higgs bundle is said to be 
 {\em semistable} if $\mu(F) \leq \mu(E)$ for any nonzero $\varphi$-invariant subbundle $F \subseteq E$. A Higgs bundle is said to be {\em stable} if $\mu(F) < \mu(E)$ for any nonzero $\varphi$-invariant subbundle $F \subsetneq E$. Finally, $(E,\varphi)$ is called {\em polystable} if it is the direct sum of stable $\varphi$-invariant subbundles, all of them of the same slope.
\end{Def}

 With this notion of stability in mind, Hitchin \cite{Hi1} constructs
the moduli space of Higgs bundles as the quotient
   $$
   \sM^{H}(2,1) = \{(E,\varphi)|\ E \text{ polystable } \} / \gG^{\CC}
   $$
and the subspace
   $$
   \sM^{H}_{s}(2,1) = 
   \{(E,\varphi)|\ E \text{ stable } \} / \gG^{\CC}
   \subseteq
   \sM^{H}(2,1),
   $$
of stable Higgs bundles. Again, given that $\GCD(2,1) = 1$, we have $\sM^{H}(2,1) = \sM^{H}_{s}(2,1)$. 
Then, Hitchin proves that the moduli of solutions of SDE coincides with the algebraic construction of the Higgs bundles moduli space. 

 \begin{thm}[\cite{Hi1}]
  There is a homeomorphism of topological spaces
  $$
   \sM^{H}(2,1)
   \cong
   \sM^{YM}(2,1).
  $$
 \end{thm}

Then, we remove the superscript and call it $ \sM (2,1)$, whose dimension was also calculated by Hitchin:

\begin{thm}\textnormal{\cite[Theorem 5.8.]{Hi1}}\label{hit5.8.}
 Let $X$ be a compact Riemann surface of genus $g > 1$. The moduli space $\sM(2,1)$ of rank $2$ and degree $1$ Higgs bundles $(E,\varphi)$, is a smooth real 
 manifold of dimension $\dim_{\RR}\sM(2,1) = 12(g - 1)$ and 
a quasi-projective variety of complex dimension
$\dim_{\CC}\sM(2,1) = 3(2g - 2)$.
\end{thm}

 \subsection[Higher rank and higher dimensional Higgs bundles]{Higher rank and dimensional Higgs bundles} 
 \label{ssec:Higgsgeneral}

Higgs bundles generalize to higher rank $r>2$. Nitsure \cite{Ni} constructs the moduli space $\mathcal{M}(r,d,L)$ of semistable pairs, which are
$$(E, \varphi: E\rightarrow E\otimes L)$$
where $E$ is a rank $r$ and degree $d$ holomorphic vector bundle over a smooth projective algebraic curve $X$, and $L$ is a line bundle over $X$. It is a GIT construction generalizing the ideas of the moduli space of vector bundles in Section \ref{ssec:GITms}. The dimension of these moduli now depend on the line bundle L, and for the case $L=K$ we obtain:
 
 \begin{thm}\textnormal{\cite{Ni}}
 \label{thm:dimHiggs}
  The space $\sM(r,d, K)$ is a quasi--projective variety of complex dimension
  \[
  \dim_{\CC}\sM(r,d, K) = (r^2-1)(2g - 2).
  \]
\end{thm}

Note that the result of Nitsure \cite{Ni} coincides with the result of Hitchin \cite{Hi1} for rank two.

Given $X$ a smooth complex projective variety of dimension $n$, a \emph{Higgs bundle} is a pair $(E,\varphi)$ where $E$ is a
locally free sheaf over $X$
and a \emph{Higgs field} $\varphi:E\rightarrow E\otimes\Omega^{1}_{X}$ verifying
$\varphi\wedge\varphi=0$. It can be thought as a coherent sheaf $\mathcal{E}$ on
the cotangent bundle $T^{\ast}X$. This definition is due to Simpson, who constructs a moduli space for  semistable Higgs bundles 
with this point of view in \cite{Si2}, what is known as the non-abelian Hodge theory.

\section{Unstability correspondence}
\label{AZresults}

In this section we survey the results of the first named author on different correspondences of the stability notion and the GIT picture, at the level of maximal unstability provided by the Harder-Narasimhan filtration.

\subsection{Correspondence for vector bundles}
\label{ssec:uns-cor}

Recall from subsection \ref{ssec:GITms} the GIT construction of moduli space trying to classify some type of
algebro-geometric objects modulo an equivalence relation. Usually,
we have to impose stability conditions on the objects we classify,
in order to obtain a space with good properties where each point
corresponds to an equivalence class of objects. By rigidifying the
objects, which typically involves adding a piece of data to the
object we want to parameterize (in the example of the moduli space of bundles, this piece of data is the isomorphism between $V$ and $H^{0}(E(m))$), we realize them as points in a
finite dimensional parameter space. The freedom in the choice of
the additional data corresponds to the action of a group.
Mumford's GIT \cite{Mu} enables then to undertake a
quotient, obtaining a projective variety which is the moduli space
classifying the objects in the moduli problem.
In every moduli problem using GIT, at some point, one has to prove
that both notions of stability do coincide, then the semistable
objects correspond to GIT semistable points, and the unstable ones
are related to the GIT unstable ones. This eventual identification, in the moduli of vector bundles, happens for a large value of the twist $m$. 

By the Hilbert-Mumford
criterion (c.f. Theorem \ref{HMcrit}), we can characterize GIT stability
through a numerical function on $1$-parameter subgroups, which turns out to be positive or
negative when the $1$-parameter subgroup destabilizes a point
or not, in the sense of GIT. Besides, when a point is GIT
unstable, we are able to talk about \emph{degrees of
unstability}, corresponding to $1$-parameter
subgroups which are more destabilizing than others. Based on the
work of Mumford, Tits and Kempf (see \cite{Ke}), among others, we
can measure this by means of a rational function on the space of
$1$-parameter subgroups, whose numerator is the numerical function
of the Hilbert-Mumford criterion and the denominator is a norm for
the $1$-parameter subgroup. We choose this norm to avoid rescaling
of the numerical function. By a theorem of Kempf \cite{Ke},
there exists a unique $1$-parameter subgroup giving a maximum for
this function, representing the maximal GIT unstability, an idea which can also be seen from the differential and symplectic point of view as the direction of maximal descense (recall subsection \ref{ssec:maximal-unstability} and the references  \cite{Ki, KN, GRS}).  Hence, an
unstable object gives a GIT unstable point for which there exists
a unique $1$-parameter subgroup GIT maximally destabilizing. From this $1$-parameter subgroup it is possible to construct a filtration by subobjects of the original unstable object,
which makes sense to ask whether it coincides with the Harder-Narasimhan filtration in cases where it is already known, or if it is able to provide a new notion of such
filtration in other cases.

In \cite{GSZ1} the main case of this correspondence is discussed, for torsion free coherent sheaves over projective
varieties, which is the generalization of the moduli of vector bundles over projective curves to higher dimensional varieties. The construction is due to Gieseker \cite{Gi1} and Maruyama \cite{Ma}, and the main difference is that the stability condition is expressed by means of a polynomial, the Hilbert polynomial $P_{E}(m)=\chi(E(m))$, encoding the information of all Chern classes, not just the degree. A torsion free coherent sheaf is said to be semistable if for every subsheaf $F\subsetneq E$, 
$$\frac{P_{F}(m)}{\rk F} \leq \frac{P_{E}(m)}{\rk E}\; .$$
Similarly to the construction in subsection \ref{ssec:GITms},  there is an isomorphism between a complex vector space $V$ and the space of global sections $H^{0}(E(m))$, yielding an action of $SL(V)$ in the Quot scheme. The moduli space is the GIT quotient by the action of this group and, during the GIT process, $1$-parameter subgroups distinguish GIT stable and GIT unstable points. A $1$-parameter
subgroup produces a weighted flag of vector subspaces 
$$0\subset V_{1}^{m}\subset V_{2}^{m}\subset \cdots V_{t}^{m}\subset V_{t+1}^{m}=V\; ,$$
which, by the isomorphism $V\simeq H^{0}(E(m))$, is a filtration of vector subspaces of global sections which can be evaluated to produce a filtration of
subsheaves of $E$, 
$$0\subset E_{1}^{m}\subset E_{2}^{m}\subset \cdots E_{t}^{m}\subset E_{t+1}^{m}=E\; .$$
 The whole process depends on an integer $m$ related to the
embedding of the Quot scheme on a projective space.
\begin{thm}\cite[Proposition 5.3 and Corollary 6.4]{GSZ1}
Given an $m$ sufficiently large, the filtration obtained by evaluating the weighted flag of vector subspaces $V_{\bullet}\subset V\simeq H^{0}(E(m))$ does not depend on $m$. Moreover, this filtration coincides with the
Harder-Narasimhan filtration of an unstable sheaf $E$ (see subsection \ref{ssec:HN}).
\end{thm}

\subsection{Other correspondences}
\label{ssec:uns-cor-others}

The previous construction can be carried out for other moduli problems. 
We call a \emph{holomorphic pair} to 
$$(E,\varphi:E\to
\sO_{X})$$
consisting on a rank $r$ vector bundle $E$ with
fixed determinant $\det(E)\cong L$ over $X$ a smooth complex
projective variety, and a morphism $\varphi$ to the trivial line
bundle $\sO_{X}$. There is a notion of stability depending
on a parameter which, in this higher dimensional case, is a polynomial. The same
techniques as in the case of sheaves apply to show that the
$1$-parameter subgroup GIT maximally destabilizing corresponds to
the Harder-Narasimhan filtration for an unstable pair (c.f. \cite{GSZ2}).
In \cite[Section 2.3]{Za1} it is shown an analogous correspondence for
the moduli problem of Higgs bundles (see Section \ref{ssec:Higgsgeneral}), following Simpson's construction \cite{Si2}.

There are moduli problems where there is no notion of what a Harder-Narasimhan filtration should be a priori. In these cases, the previous correspondence can help to define a Harder-Narasimhan filtration as that filtration coming from the maximal destabilizing $1$-parameter subgroup, from the GIT point of view. 
We call a \emph{rank $2$ tensor} the pair
consisting of 
$$
(E,\; \varphi:{E\otimes\cdots \otimes E}\too M)\; ,$$ where
$E$ is a rank $2$ coherent torsion free sheaf over a smooth
complex projective variety $X$ and $M$ a line bundle over $X$. In this case, the Harder-Narasimhan filtrations are simply line subbundles $L$.
In \cite{Za3}, symmetric rank $2$ tensors are interpreted as degree $s$
coverings $X'\rightarrow X$ lying on the ruled surface
$\mathbb{P}(E)$, to define a notion of stable covering and
characterize geometrically the maximally destabilizing subbundle
$L\subset E$ in terms of intersection theory and configurations of points as in Example \ref{npoints}.

There are also other categories, such as quiver representations, where these ideas apply. 
Let $Q$ be a finite quiver and consider a representation of $Q$ on finite
dimensional $k$-vector spaces, where $k$ is an algebraically
closed field of arbitrary characteristic. We consider the construction of a moduli space for these objects
by King \cite{Ki} and associate to an unstable
representation an unstable point, in the sense of GIT,
in a parameter space where a group acts. Then,
the $1$-parameter subgroup gives a filtration of
subrepresentations and we prove that it coincides with the
Harder-Narasimhan filtration for that representation (see \cite{Za3} for quiver representations and \cite[Chapter 3]{Za1} for Q-sheaves).

Many of the cases where we have been able to carry this
correspondence out fall into abelian categories, where the existence of 
Harder-Narasimhan filtrations is straightforward. In other cases where the moduli construction is performed in a non-abelian category, 
existence and uniqueness of a Harder-Narasimhan filtration is a much harder problem. In \cite{TZ} it is studied the moduli problem of $(G,h)$-constellations, whose stability condition can be seen as an infinite dimensional version of the one for quiver representations. There, it is shown that the maximal unstable filtration coming from the GIT picture does not necessarily stabilize with the embedding (depending on a subset $D$ of the set of irreducible representations of $G$) in certain Quot scheme. However, it produces a filtration that asymptotically converges to the Harder-Narasimhan filtration, which is proved to exist. 

Coming from these ideas it could be thought that a similar construction can be performed for principal $G$-bundles. The notion of Harder-Narasimhan filtration for principal bundles is achieved by means of the canonical reduction, a reduction to a parabolic subgroup $P\subset G$. In \cite{BZ},
examples of orthogonal and symplectic bundles over surfaces are shown where the
Harder-Narasimhan filtration of its underlying vector bundle does not correspond
to any parabolic reduction of the orthogonal or symplectic bundle. This result prevents us from trying to relate unstability of $G$-bundles with unstability of the underlying bundle so straightforwardly, as it happens with stability.

\section[Stratifications on the moduli space of Higgs bundles]{Stratifications on the moduli space of Higgs bundles}
\label{sec:RZRresults}

In this final section we review results of the second named author on stratifications of the moduli space of Higgs bundles, performed with the invariants provided by their Harder-Narasimhan filtrations.

\subsection{Shatz stratification}
\label{ssec:shatz}

\begin{figure}[h!]
\begin{center}
  \includegraphics[scale=0.5]{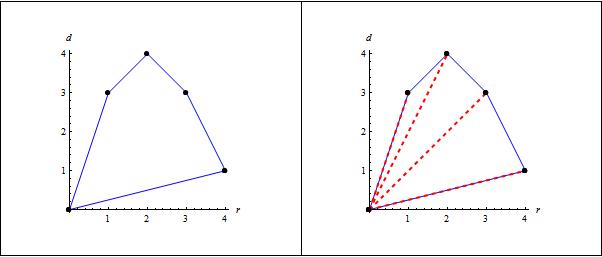}
\caption{Harder-Narasimhan poligons.}
\label{HNpolygons}
\end{center}
\end{figure}

Recall the Harder-Narasimhan filtration of a holomorphic vector bundle, from Section \ref{ssec:HN}. Any vector bundle has a unique Harder-Narasimhan filtration:
$$ \HNF(E): \left\{ 0\right\}=E_{0}\subsetneq E_{1}\subsetneq E_{2}\subsetneq ...\subsetneq E_{t-1}\subsetneq E_{t}=E\; .$$
We define the \emph{Harder-Narasimhan polygon} 
as the polygon in the first quadrant of the $(r,d)$-plane with vertices $(\rk(E_i),\deg(E_i))$ for $i=0,\dots,t$ (see the black points in Figure \ref{HNpolygons}, left). 
The slope of the line joining $(\rk(E_{i-1}),\deg(E_{i-1}))$ and $(\rk(E_i),\deg(E_i))$ is $\mu(E_i/E_{i-1})=\mu(E^i)$ (blue line segments in Figure \ref{HNpolygons}, left).  
Slopes $\mu(E_i)$ of each factor in the filtration are represented by red-dashed lines in Figure \ref{HNpolygons}, right. The \emph{Harder-Narasimhan type} of $E$ is defined as the vector in $\QQ^r$:
$$
  \HNT(E) = \overrightarrow{\mu} = \big(\mu({E^1}),\dots,\mu({E^1}),\dots,\mu(E^t),\dots,\mu(E^t)\big)
$$
where the slope of each $E^j$ appears $\rk(E^j)$ times.

First condition in Theorem \ref{HNdef} states that $\mu(E_{i}) < \mu(E_{i-1})$ for $i=2,\dots,s$, being equivalent to say that the Harder-Narasimhan polygon is convex, which is intuitively clear from the definiton of the Harder-Narasimhan polygon (see \cite[Proposition 5]{Sha}). Second condition assures that all quotients are semistable, which can be rephrased as the Harder-Narasimhan polygon being the convex envelope of every refinement (i.e. if there is an unstable quotient, it has a destabilizing subbundle of greater slope which therefore will be depicted above the polygon).

Let $\sM(r,d)$ be the moduli space of Higgs bundles $(E, \varphi:E\rightarrow E\otimes K)$, where $E$ has rank $r$ and degree $d$. 
 As a consequence of the work of Shatz \cite[Propositions 10~\&~11]{Sha}, there is a finite stratification of $\sM(r,d)$ by the Harder-Narasimhan type of the underlying vector bundle $E$:
 $$
 \sM(r,d) = \bigcup_{\overrightarrow{\mu}}U'_{\overrightarrow{\mu}}
 $$ 
 where $U'_{\overrightarrow{\mu}} \sse \sM(r,d)$ is the subspace of Higgs bundles $(E,\varphi)$ whose underlying vector bundle $E$ has Harder-Narasimhan type $\overrightarrow{\mu}$. This stratification is known as the {\em Shatz stratification}. Note that there is an open dense stratum corresponding to Higgs bundles $(E,\varphi)$ such that the underlying vector bundle $E$ is itself stable. Since $\varphi\in H^0(\End(E)\otimes K)\cong H^1(\End(E))^*$ (by Serre duality) such Higgs bundles represents a point in the cotangent bundle of the moduli space of stable bundles $\mathcal{N}(r,d)$, namely
$$
  U'_{(d/r,\dots,d/r)} = T^*\mathcal{N}(r,d)\sse \sM(r,d)\; .
$$

\subsection[$\mathbb{C}^{\ast}$-action and Bia\l{}ynicki-Birula stratification]{$\mathbb{C}^{\ast}$-action and Bia\l{}ynicki-Birula stratification}
\label{ssec:BBstrat}

The moduli space of Higgs bundles $\sM(r,d)$ is endowed with an important action of the group $\mathbb{C}^{\ast}$, defined by multiplication (see \cite{Si2})
$$
z\cdot (E,\varphi)\mapsto (E,z\cdot \varphi).
$$
The limit $\displaystyle (E_0,\varphi_0)=\lim_{z\to0}z\cdot (E, \varphi)$ exists for all $(E,\varphi)\in
\sM(r,d)$ and, moreover, this limit is a fixed boint by the $\CC^{*}$-action. A Higgs bundle $(E,\varphi)$ is a fixed 
point of the $\CC^*$-action if and only if it is a \emph{Hodge bundle}: there is a decomposition
$$
  E=\bigoplus_{j=1}^p E_j
$$
with respect to which the Higgs field $\varphi$ has weight one, meaning that $\varphi\colon E_j\to E_{j+1}\otimes K$, for all $j$. The 
\emph{type} of the Hodge bundle is the sequence of ranks of $(E,\varphi)$, $\big(\rk(E_1),\dots,\rk(E_p)\big)$. Let $\{F_{\lambda}\}$ be 
the irreducible components of the fixed point locus of $\CC^{*}$ on $\sM(r,d)$.  Consider the upper flow sets
$$
 U_{\lambda}^{+} := 
\{ (E, \varphi)\in \sM \: \lim_{z\to0}z\cdot (E, \varphi) \in F_{\lambda} \}.
$$
Then we have the {\em Bia\l{}ynicki-Birula stratification}  of $\sM(r,d)$   (cf. \cite{B-B}),
$$
 \sM = \bigcup_{\lambda}U_{\lambda}^+\; ,
$$
which is indexed by these irreducible components. 

\subsection[Moduli space stratifications]{Moduli space stratifications}
\label{ssec:Higgsstrat}

Motivated by the results of Hausel \cite{Hau} on stratifications of the moduli space of rank two Higgs bundles, there has been some developments in stratifications for rank three (detailed proofs can be found in  \cite{GZ-R}). 

First, we 
present bounds on Harder-Narasimhan types for rank three Higgs bundles. Let $(E,\varphi)$ be a rank $3$ Higgs bundle, and let $(\mu_1,\mu_2,\mu_3)$ be the Harder-Narasimhan type of the underlying bundle $E$, so that $\mu_1\geq\mu_2\geq\mu_3$ and $\mu_1+\mu_2+\mu_3=3\mu$, where $\mu=\mu(E)$. We can write the Harder-Narasimhan filtration of the vector bundle $E$ as follows:
$$
  \HNF(E): \{0\} = E_0 \sse E_1 \sse E_2 \sse E_3 = E,
$$
where $E_i=E_j$ if
$\mu_i=\mu_j$, by convention. Thus, for example, if $\mu_1=\mu_2>\mu_3$ then the
Harder--Narasimhan filtration is
$$
  \HNF(E): \{0\}= E_0 \sse E_1 = E_2 \sse E_3 = E
$$
and $\rk(E_1)=\rk(E_2)=2$. Similarly, if $\mu_1>\mu_2=\mu_3$, then
$\rk(E_1)=1$ and $\rk(E_2)=3$.

\begin{rem}
  \label{rmk:phi-non-zero}
  Let $(E,\varphi)$ be a stable Higgs bundle such that $E$ is an unstable
  vector bundle of $\HNT(E)=(\mu_1,\mu_2,\mu_3)$. Then
  $E_1\sse E_2$ is destabilizing and hence, by stability of
  $(E,\varphi)$, we have $\varphi_{21}\neq 0$ where 
  $
  \varphi = 
  \left(
    \begin{smallmatrix}
       0      &       0      & 0\\
         \varphi_{21} &       0      & 0\\
               0      & \varphi_{32} & 0
    \end{smallmatrix}
  \right)$. 
  Similarly $E_2\sse E$
  is destabilizing and therefore $\varphi_{32}\neq 0$ (unless
  $\mu_2=\mu_3 \iff E_2=E$).
\end{rem}

 After analyzing the image and the kernel of $\varphi$, we get the following bounds:

\begin{prop}
\textnormal{\cite[Proposition 4.2]{GZ-R}}
  \label{prop:mu-inequalities}
  Let $(E,\varphi)$ be a semistable rank $3$ Higgs bundle of
  $\HNT(E)=(\mu_1,\mu_2,\mu_3)$. Then
  \begin{align}
    0 &\leq \mu_1-\mu_2 \leq 2g-2, \label{eq:mu1-mu2}\\
    0 &\leq \mu_2-\mu_3 \leq 2g-2. \label{eq:mu2-mu3}
  \end{align}
\end{prop}

 The purpose now will be to analyze the limit of the $\CC^{*}$-action as $z\to 0$ in terms of $\HNT(E)$. 
Let us start by trivial filtrations. Let $(E,\varphi)$ be a stable Higgs bundle. When the underlying vector
bundle $E$ is also itself stable, clearly $\displaystyle \lim_{z\to 0}z\cdot(E,\varphi)=(E,0)$. Hence we have the following result, valid
for any rank.

\begin{prop}
\textnormal{\cite[Proposition 2.2]{GZ-R}}
  Let $(E,\varphi)\in\mathcal{M}(r,d)$. Then $\displaystyle \lim_{z\to 0}z\cdot(E,\varphi)=(E,0)$ if and only if $E$ is stable.
\end{prop}

 Now, let us check on non-trivial Harder-Narasimhan filtrations. Consider rank $3$
stable Higgs bundles $(E,\varphi)$ and assume that $GCD(3,d)=1$. 

\begin{thm}
\textnormal{\cite[Theorem 5.1]{GZ-R}}
\label{stratifications}
  Let $(E,\varphi) \in \mathcal{M}(3,d)$ be such that $E$ is an unstable
  vector bundle of slope $\mu$ and 
  $\HNT(E)=(\mu_1,\mu_2,\mu_3)$. Then the limit 
  $
  \displaystyle
  (E_0,\varphi_0)=
  \lim_{z\to0}(E,z\cdot \varphi)
  $ 
  falls into one of these cases:

  \begin{enumerate}
  \item[(1)] Assume that $\mu_2<\mu$. Then $\mu_1>\mu_2\geq \mu_3$
    and one of the following alternatives holds.
    \begin{enumerate}
    \item[(1.1)] The slope of $I=\varphi_{21}(E_1)\ox K^{-1}$ satisfies $\mu_1-(2g-2) \leq
      \mu(I)<-\frac{1}{3}\mu_1+\frac{2}{3}\mu_2+\frac{2}{3}\mu_3$ and
      $(E_0,\varphi_0)$ is the following Hodge bundle of type $(1,2)$:
      $$
      (E_0,\varphi_0) = \Big(E_1\oplus E/E_1,
      \left(
       \begin{smallmatrix}
               0      &       0\\
         \varphi_{21} &  0
       \end{smallmatrix}
      \right) 
      \Big).
      $$
      The associated graded vector bundle is $\Gr(E_0) = \Gr(E)$.
    \item[(1.2)] The slope of $I=\varphi_{21}(E_1)\ox K^{-1}$ satisfies 
      $-\frac{1}{3}\mu_1+\frac{2}{3}\mu_2+\frac{2}{3}\mu_3<\mu(I)\leq
      \mu_3$ and $(E_0,\varphi_0)$ is the following Hodge bundle of type
      $(1,1,1)$:
      $$
      (E_0,\varphi_0) = \Big(
      E_1\oplus I \oplus (E/E_1)/I, 
      \left(
      \begin{smallmatrix}
               0      &       0      & 0\\
         \varphi_{21} &       0      & 0\\
               0      & \varphi_{32} & 0
      \end{smallmatrix}
     \right) 
     \Big),
     $$
     where $\varphi_{21}$ and $\varphi_{32}$ are induced from $\varphi$.
     The associated graded vector bundle is
     $\Gr(E_0) = E_1 \oplus (E/E_1)/I \oplus I$ and 
     $\HNT(E_0)=(\mu_1,\mu_2+\mu_3-\mu(I),\mu(I))$.
   \item[(1.3)] The slope of $I=\varphi_{21}(E_1)\ox K^{-1}$ satisfies $\mu(I)=\mu_2$ and the
     strict inequality $\mu_3<\mu_2$ holds. Moreover, the line bundle
     $I=E_2/E_1$ and $(E_0,\varphi_0)$ is the following Hodge bundle of
     type $(1,1,1)$:
      $$
      (E_0,\varphi_0) = 
      \Big(
      E_1\oplus E_2/E_1 \oplus E/E_2, 
      \left(
       \begin{smallmatrix}
               0      &       0      & 0\\
         \varphi_{21} &       0      & 0\\
               0      & \varphi_{32} & 0
       \end{smallmatrix}
      \right) 
      \Big),
     $$ 
     where $\varphi_{32}$ is induced from $\varphi$.
     The associated graded vector bundle is $\Gr(E_0) = \Gr(E)$.
    \end{enumerate}
  \item[(2)] Suppose that $\mu_2>\mu$. Then $\mu_1\geq\mu_2>\mu_3$ and
    one of the following alternatives holds.
    \begin{enumerate}
    \item[(2.1)] The slope of $N = \ker(\varphi_{32})$ satisfies $\mu_1+\mu_2-\mu_3-(2g-2)
\leq \mu(N)<\mu$ and $(E_0,\varphi_0)$ is the following Hodge
      bundle of type $(2,1)$:
      $$
      (E_0,\varphi_0) = \Big(E_2\oplus E/E_2,
      \left(
    \begin{smallmatrix}
               0      & 0\\
         \varphi_{32} & 0
    \end{smallmatrix}
      \right) 
      \Big).
      $$
      The associated graded vector bundle is $\Gr(E_0) = \Gr(E)$.
    \item[(2.2)] The slope of $N = \ker(\varphi_{32})$ satisfies
      $\mu<\mu(N)\leq\mu_2$ and $(E_0,\varphi_0)$ is the following
      Hodge bundle of type $(1,1,1)$:
      $$
      (E_0,\varphi_0) = \Big(N \oplus {E}_2/N \oplus E/E_2,
      \left(
       \begin{smallmatrix}
               0      &       0      & 0\\
         \varphi_{21} &       0      & 0\\
               0      & \varphi_{32} & 0
       \end{smallmatrix}
      \right) 
                     \Big)
      $$
      where $\varphi_{21}$ and $\varphi_{32}$ are induced from
      $\varphi$. The associated graded vector bundle is
      $\Gr(E_0) = E_2/N \oplus N \oplus E/E_2$ and 
      $\HNT(E_0) = (\mu_1+\mu_2-\mu(N),\mu(N),\mu_3)$.
    \item[(2.3)] The slope of $N = \ker(\varphi_{32})$ satisfies $\mu(N)=\mu_1$ and the strict
      inequality $\mu_1>\mu_2$ holds. Moreover the line bundle $N=E_1$
      and $(E_0,\varphi_0)$ is the following Hodge bundle of type
      $(1,1,1)$:
      $$
      (E_0,\varphi_0) = 
      \Big(
      E_1\oplus E_2/E_1 \oplus E/E_2, 
      \left(
       \begin{smallmatrix}
               0      &       0      & 0\\
         \varphi_{21} &       0      & 0\\
               0      & \varphi_{32} & 0
       \end{smallmatrix}
      \right) 
      \Big),
     $$ 
     where $\varphi_{21}$ and $\varphi_{32}$ are induced from $\varphi$.
     The associated graded vector bundle is $\Gr(E_0) = \Gr(E)$.
    \end{enumerate}
  \end{enumerate}
\end{thm}

\begin{cor}
\textnormal{\cite[Corollary 5.4]{GZ-R}}
\label{cor:111a}
  Let $(E,\varphi) \in \mathcal{M}(3,d)$ be such that $E$ is an unstable
  vector bundle of slope $\mu$ and 
  $\HNT(E)=(\mu_1,\mu_2,\mu_3)$. Assume that $\mu_1-\mu_3>2g-2$. Then the
  limit $\displaystyle (E_0,\varphi_0)= \lim_{z\to0}(E,z\cdot \varphi)$ is given by
  (1.3) of Theorem~\ref{stratifications} if $\mu_2<\mu$, and by
  (2.3) of Theorem~\ref{stratifications} if $\mu_2>\mu$.
\end{cor}

\begin{cor}
\textnormal{\cite[Corollary 5.5]{GZ-R}}
\label{cor:111b}
  Let $(E_0=L_1\oplus L_2\oplus L_3,\varphi_0=\left(
    \begin{smallmatrix}
       0      &       0      & 0\\
         \varphi_{21} &       0      & 0\\
               0      & \varphi_{32} & 0
    \end{smallmatrix}
\right))$ be a Hodge bundle
  of type $(1,1,1)$ with $\mu(L_1)-\mu(L_3)>2g-2$. Then
  $\mu(L_1)>\mu(L_2)>\mu(L_3)$ and any $(E,\varphi)$ such that
  $\displaystyle \lim_{z\to0}(E,z\cdot \varphi) = (E_0,\varphi_0)$ satisfies
  $\Gr(E)=\Gr(E_0)$.
\end{cor}

Corollaries \ref{cor:111a} and \ref{cor:111b} lead to an identification between Shatz strata and Bia\l{}ynicki-Birula strata in some cases. Recall that the set $U^+_{(l_1,l_2,l_3)}$ denotes the Bia\l{}ynicki-Birula stratum of Higgs bundles whose limits lie in $F_{(l_1,l_2,l_3)}$ and that $U'_{(l_1,l_2,l_3)}$ denotes the Shatz stratum of Higgs bundles whose Harder--Narasimhan type is $(l_1,l_2,l_3)$.

\begin{thm}
\textnormal{\cite[Theorem 5.6]{GZ-R}}
  \label{thm:HN=BB}
  Let $(l_1,l_2,l_3)$ be such that $l_1-l_3>2g-2$. Then the
  corresponding Shatz and Bia\l{}ynicki-Birula strata in
  $\mathcal{M}(3,d)$ coincide:
  \begin{displaymath}
    U'_{(l_1,l_2,l_3)} = U^+_{(l_1,l_2,l_3)}.
  \end{displaymath}
\end{thm}

\end{document}